\documentclass[aos,nameyear,dvips]{arximspdf}
\usepackage{dcolumn}
\usepackage{graphicx}

\doi{10.1214/09-AOS762}
\volume{38}
\issue{4}
\pubyear{2010}
\firstpage{2005}
\lastpage{2046}

\makeatletter

\newcolumntype{d}[1]{D{.}{.}{#1}}

\newtheorem{lemma}{Lemma}

\newproclaim{remark}{Remark}

\newtheorem{theorem}{Theorem}
\newtheorem{proposition}{Proposition}

\makeatother

\begin{document}
\begin{frontmatter}

\title{Nonparametric regression in exponential families}
\runtitle{Nonparametric regression in exponential families}

\begin{aug}
\author[A]{\fnms{Lawrence D.} \snm{Brown}\thanksref{t1}},
\author[A]{\fnms{T. Tony} \snm{Cai}\corref{}\thanksref{t2}\ead[label=e1]{tcai@wharton.upenn.edu}} and
\author[B]{\fnms{Harrison H.} \snm{Zhou}\thanksref{t3}\ead[label=e3]{huibin.zhou@yale.edu}}
\runauthor{L. D. Brown, T. T. Cai and H. H. Zhou}
\affiliation{University of Pennsylvania, University of Pennsylvania
and Yale University}
\address[A]{L. D. Brown\\
T. T. Cai\\
Department of Statistics\\
The Wharton School\\
University of Pennsylvania\\
Philadelphia, Pennsylvania 19104\\
USA\\
\printead{e1}}
\address[B]{H. H. Zhou\\
Department of Statistics\\
Yale University\\
P.O. Box 208290\\
New Haven, Connecticut 06520-8290\\
USA\\
\printead{e3}}
\end{aug}

\thankstext{t1}{Supported in part by NSF Grant DMS-07-07033.}

\thankstext{t2}{Supported in part by NSF Grant DMS-06-04954
and NSF FRG Grant DMS-08-54973.}

\thankstext{t3}{Supported in part by NSF Career Award DMS-06-45676 and
NSF FRG Grant DMS-08-54975.}

\received{\smonth{2} \syear{2009}}
\revised{\smonth{10} \syear{2009}}

%
\begin{abstract}
Most results in nonparametric regression theory are developed only for the
case of additive noise. In such a setting many smoothing techniques
including wavelet thresholding methods have been developed and shown to be
highly adaptive. In this paper we consider nonparametric regression in
exponential families with the main focus on the natural exponential families
with a quadratic variance function, which include, for example, Poisson
regression, binomial regression and gamma regression. We propose a unified
approach of using a mean-matching variance stabilizing transformation to
turn the relatively complicated problem of nonparametric regression in
exponential families into a standard homoscedastic Gaussian regression
problem. Then in principle any good nonparametric Gaussian regression
procedure can be applied to the transformed data. To illustrate our general
methodology, in this paper we use wavelet block thresholding to construct
the final estimators of the regression function. The procedures are easily
implementable. Both theoretical and numerical properties of the estimators
are investigated. The estimators are shown to enjoy a high degree of
adaptivity and spatial adaptivity with near-optimal asymptotic performance
over a wide range of Besov spaces. The estimators also perform well
numerically.
\end{abstract}

%

\end{frontmatter}

\section{Introduction}
\label{intro}

Theory and methodology for nonparametric regression is now well developed
for the case of additive noise particularly additive homoscedastic Gaussian
noise. In such a setting many smoothing techniques including wavelet
thresholding methods have been developed and shown to be adaptive and enjoy
other desirable properties over a wide range of function spaces.
However, in
many applications the noise is not additive and the conventional
methods are
not readily applicable. For example, such is the case when the data are
counts or proportions.

In this paper we consider nonparametric regression in exponential families
with the main focus on the natural exponential families with a quadratic
variance function (NEF--QVF). These include, for example, Poisson regression,
binomial regression and gamma regression. We present a unified
treatment of
these regression problems by using a mean-matching variance stabilizing
transformation (VST) approach. The mean-matching VST turns the relatively
complicated problem of regression in exponential families into a standard
homoscedastic Gaussian regression problem and then any good nonparametric
Gaussian regression procedure can be applied.

Variance stabilizing transformations and closely related normalizing
transformations have been widely used in many parametric statistical
inference problems. See Hoyle (\citeyear{Hoyle73}), Efron (\citeyear{Efron82}) and Bar-Lev and Enis
(\citeyear{BE90}). In the more standard parametric problems, the goal of VST is often
to optimally stabilize the variance. That is, one desires the variance of
the transformed variable to be as close to a constant as possible. For
example, Anscombe (\citeyear{Anscombe48}) introduced VSTs for binomial, Poisson and negative
binomial distributions that provide the greatest asymptotic control
over the
variance of the resulting transformed variables. In the context of
nonparametric function estimation, Anscombe's variance stabilizing
transformation has also been briefly discussed in Donoho (\citeyear{Donoho93}) for density
estimation. However, for our purposes it is much more essential to have
optimal asymptotic control over the bias of the transformed variables. A
mean-matching VST minimizes the bias of the transformed data while also
stabilizing the variance.

Our procedure begins by grouping the data into many small size bins,
and by
then applying the mean-matching VST to the binned data. In principle any
good Gaussian regression procedure could be applied to the transformed data
to construct the final estimator of the regression function. To illustrate
our general methodology, in this paper we employ two wavelet block
thresholding procedures. Wavelet thresholding methods have achieved
considerable success in nonparametric regression in terms of spatial
adaptivity and asymptotic optimality. In particular, block thresholding
rules have been shown to possess impressive properties. In the context of
nonparametric regression, local block thresholding has been studied, for
example, in Hall, Kerkyacharian and Picard (\citeyear{HKP98}), Cai (\citeyear{Cai99}, \citeyear{Cai02}) and Cai
and Silverman (\citeyear{CS01}). In this paper we shall use the BlockJS procedure
proposed in Cai (\citeyear{Cai99}) and the NeighCoeff procedure introduced in Cai and
Silverman (\citeyear{CS01}). Both estimators were originally developed for
nonparametric Gaussian regression. BlockJS first divides the empirical
coefficients at each resolution level into nonoverlapping blocks and then
simultaneously estimates all the coefficients within a block by a
James--Stein rule. NeighCoeff also thresholds the empirical
coefficients in
blocks, but estimates wavelet coefficients individually. It chooses a
threshold for each coefficient by referencing not only to that coefficient
but also to its neighbors. Both estimators increase estimation accuracy over
term-by-term thresholding by utilizing information about neighboring
coefficients.

Both theoretical and numerical properties of our estimators are
investigated. It is shown that the estimators enjoy excellent asymptotic
adaptivity and spatial adaptivity. The procedure using BlockJS
simultaneously attains the optimal rate of convergence under the integrated
squared error over a wide range of the Besov classes. The estimators also
automatically adapt to the local smoothness of the underlying function; they
attain the local adaptive minimax rate for estimating functions at a point.
A key step in the technical argument is the use of the quantile coupling
inequality of Koml\'{o}s, Major and Tusn\'{a}dy (\citeyear{KMT75}) to approximate the
binned and transformed data by independent normal variables. The procedures
are easy to implement, at the computational cost of $O(n)$. In addition to
enjoying the desirable theoretical properties, the procedures also perform
well numerically. 

Our method is applicable in more general settings. It can be extended to
treat nonparametric regression in general one-parameter natural exponential
families. The mean-matching VST only exists in NEF--QVF (see Section
\ref%
{vst.sec}). In the general case when the variance is not a quadratic
function of the mean, we apply the same procedure with the standard VST in
place of the mean-matching VST. It is shown that, under slightly stronger
conditions, the same optimality results hold in general. We also note that
mean-matching VST transformations exist for some useful
nonexponential families, including some commonly used for modeling
``over-dispersed'' data. Though we do not pursue the details in the
present paper, it appears that because of this our methods can also be
effectively used for nonparametric regressions involving such error
distributions.

We should note that nonparametric regression in exponential families has
been considered in the literature. Among individual exponential families,
the Poisson case is perhaps the most studied. Besbeas, De Feis and Sapatinas
(\citeyear{BFS04}) provided a review of the literature on the nonparametric Poisson
regression and carried out an extensive numerical comparison of several
estimation procedures including Donoho (\citeyear{Donoho93}), Kolaczyk (\citeyear{Kolaczyk99a}, \citeyear{Kolaczyk99b}) and
Fry\'{z}lewicz and Nason (\citeyear{FN01}). In the case of Bernoulli regression,
Antoniadis and Leblanc (\citeyear{AntLeb00}) introduced a wavelet procedure based on
diagonal linear shrinkers. Unified treatments for nonparametric regression
in exponential families have also been proposed. Antoniadis and Sapatinas
(\citeyear{AntSap01}) introduced a wavelet shrinkage and modulation method for regression
in NEF--QVF and showed that the estimator attains the optimal rate over the
classical Sobolev spaces. Kolaczyk and Nowak (\citeyear{KN05}) proposed a recursive
partition and complexity-penalized likelihood method. The estimator was
shown to be within a logarithmic factor of the minimax rate under squared
Hellinger loss over Besov spaces.

The paper is organized as follows. Section \ref{vst.sec} discusses the
mean-matching variance stabilizing transformation for natural exponential
families. In Section \ref{procedure.sec}, we first introduce the general
approach of using the mean-matching VST to convert nonparametric regression
in exponential families into a nonparametric Gaussian regression problem,
and then present in detail specific estimation procedures based on the
mean-matching VST and wavelet block thresholding. Theoretical
properties of
the procedures are treated in Section \ref{theory.sec}. Section \ref%
{numerical.sec} investigates the numerical performance of the
estimators. We
also illustrate our estimation procedures in the analysis of two real data
sets: a gamma-ray burst data set and a packet loss data set. Technical
proofs are given in Section \ref{proof.sec}.

\section{Mean-matching variance stabilizing transformation}
\label{vst.sec}

We begin by considering variance stabilizing transformations (VST) for
natural exponential families. As mentioned in the \hyperref
[intro]{Introduction}, VST has been
widely used in many contexts and the conventional goal of VST is to
optimally stabilize the variance. See, for example, Anscombe (\citeyear{Anscombe48}) and
Hoyle (\citeyear{Hoyle73}). For our purpose of nonparametric regression in exponential
families, we shall first develop a new class of VSTs, called mean-matching
VSTs, which asymptotically minimize the bias of the transformed variables
while at the same time stabilizing the variance.

Let $X_{1},X_{2},\ldots,X_{m}$ be a random sample from a distribution
in a
natural one-parameter exponential families with the probability density/mass
function
\[
q(x|\eta)=e^{\eta x-\psi(\eta)}h(x).
\]
Here $\eta$ is called the natural parameter. The mean and variance are,
respectively,
\[
\mu(\eta)=\psi^{\prime}(\eta) \quad\mbox{and}\quad \sigma^{2}(\eta)=\psi
^{\prime\prime}(\eta).
\]
We shall denote the distribution by $\operatorname{NEF}(\mu)$. A special subclass of
interest is the one with a quadratic variance function (QVF),
%
\begin{equation} \label{Q.var}
\sigma^{2}\equiv V(\mu)=a_{0}+a_{1}\mu+a_{2}\mu^{2}.
\end{equation}
In this case we shall write $X_{i}\sim \operatorname{NQ}(\mu)$. The NEF--QVF families
consist of six distributions, three continuous: normal, gamma and NEF--GHS
distributions and three discrete: binomial, negative binomial and Poisson.
See, for example, Morris (\citeyear{Morris82}) and Brown (\citeyear{Brown86}).

Set $X=\sum_{i=1}^{m}X_{i}$. According to the central limit theorem,
\[
\sqrt{m}\bigl(X/m-\mu(\eta)\bigr)\stackrel{L}{\longrightarrow} N(0,V(\mu
(\eta
))) \qquad\mbox{as }m\rightarrow\infty.
\]
A variance stabilizing transformation (VST) is a function $G\dvtx\mathbb{R}
\rightarrow\mathbb{R}$ such that
%
\begin{equation}\label{vst}
G^{\prime}(\mu)=V^{-{{1/2}}}(\mu).
\end{equation}
The standard delta method then yields
\[
\sqrt{m}\{G(X/m)-G(\mu(\eta))\}\stackrel{L}{\longrightarrow}N(0,1).
\]
It is known that the variance stabilizing properties can often be further
improved by using a transformation of the form
%
\begin{equation} \label{Hm}
H_{m}(X)=G\biggl({\frac{X+a}{m+b}}\biggr)
\end{equation}
with suitable choice of constants $a$ and $b$. See, for example,
Anscombe (\citeyear{Anscombe48}).
In this paper we shall use the VST as a tool for nonparametric
regression in
exponential families. For this purpose, it is more important to optimally
match the means than to optimally stabilize the variance. That is, we wish
to choose the constants $a$ and $b$ such that $\mathbb{E}\{H_{m}(X)\}$
optimally matches $G(\mu(\eta))$.

To derive the optimal choice of $a$ and $b$, we need the following
expansions for the mean and variance of the transformed variable $H_{m}(X)$.
\begin{lemma}
\label{approx.lem} Let $\Theta$ be a compact set in the interior of the
natural parameter space. Then for $\eta\in\Theta$ and for constants $a$
and $b$,
%
\begin{equation} \label{bias}\quad
\mathbb{E}\{H_{m}(X)\}-G(\mu(\eta))={\frac{1}{\sigma(\eta
)}}\biggl(a-b\mu
(\eta)-\frac{\mu^{\prime\prime}(\eta)}{4\mu^{\prime}(\eta
)}\biggr)\cdot
m^{-1}+O(m^{-2})
\end{equation}
and%
%
\begin{equation} \label{var}
\operatorname{Var}\{H_{m}(X)\}=\frac{1}{m}+O(m^{-2}).
\end{equation}
Moreover, there exist constants $a$ and $b$ such that
%
\begin{equation} \label{sharpbias}
\mathbb{E}\biggl\{G\biggl(\frac{X+a}{m+b}\biggr)\biggr\}-G(\mu(\eta))=O(m^{-2})
\end{equation}
for all $\eta\in\Theta$ with a positive Lebesgue measure if and only if
the exponential family has a quadratic variance function.
\end{lemma}

The proof of Lemma \ref{approx.lem} is given in Section \ref
{proof.sec}. The
last part of Lemma \ref{approx.lem} can be easily explained as follows.
Equation (\ref{bias}) implies that (\ref{sharpbias}) holds
if and
only if
\[
a-b\mu(\eta)-\frac{\mu^{\prime\prime}(\eta)}{4\mu^{\prime
}(\eta)}=0,
\]
that is, $\mu^{\prime\prime}(\eta)=4a\mu^{\prime}(\eta)-4b\mu
(\eta)\mu
^{\prime}(\eta)$. Solving this differential equation yields
%
\begin{equation} \label{sol}
\sigma^{2}(\eta)=\mu^{\prime}(\eta)=a_{0}+4a\mu(\eta)-2b\mu
^{2}(\eta)
\end{equation}
for some constant $a_0$. Hence the solution of the differential
equation is
exactly the subclass of natural exponential family with a quadratic variance
function (QVF).

It follows from (\ref{sol}) that among the VSTs of the form
(\ref%
{Hm}) for the exponential family with a quadratic variance function
\[
\sigma^{2}=a_{0}+a_{1}\mu+a_{2}\mu^{2},
\]
the best constants $a$ and $b$ for mean-matching are
%
\begin{equation}\label{optimalconstants}
a=\tfrac{1}{4}a_{1} \quad\mbox{and}\quad b=-\tfrac{1}{2}a_{2}.
\end{equation}
We shall call the VST (\ref{Hm}) with the constants $a$ and $b$ given
in (%
\ref{optimalconstants}) the mean-matching VST. Lemma \ref{approx.lem} shows
that the mean-matching VST only exists in the NEF--QVF families and
with the
mean-matching VST the bias $\mathbb{E}\{G(\frac{X+a}{m+b})\}-G(\mu
(\eta))$
is of the order $(m^{-2})$. In contrast, for an NEF without a quadratic
variance function, the term $a-\mu(\eta)b-\frac{\mu^{\prime\prime
} ( \eta) }{4\mu^{\prime} ( \eta) }$ does not vanish
for all $\eta$ with any choice of $a$ and $b$. And in this case the bias
\[
\mathbb{E}\biggl\{G\biggl(\frac{X+a}{m+b}\biggr)\biggr\}-G(\mu(\eta))=O(m^{-1})
\]
instead of $O(m^{-2})$ in (\ref{sharpbias}). We shall see in
Section \ref{theory.sec} that this difference has important
implications for
nonparametric regression in NEF.

The following are the specific expressions of the mean-matching VST $H_m$
for the five distributions (other than normal) in the NEF--QVF families:

\begin{itemize}
\item Poisson: $a=1/4$, $b=0$ and $H_{m}(X)=2\sqrt{(X+{\frac{1}{4}})/m}$.

\item $\operatorname{Binomial}( r,p ) $: $a=1/4$, $b=\frac{1}{2r}$ and $%
H_{m}(X)=2\sqrt{r}\arcsin( \sqrt{\frac{X+1/4}{rm+1/2}} ) $.

\item Negative $\operatorname{Binomial}( r,p ) $: $a=1/4$, $b=-\frac{1}{2r}$ and
\[
H_{m}(X)=2\sqrt{r}\ln\Biggl( \sqrt{\frac{X+1/4}{mr-1/2}}+\sqrt{1+{\frac
{X+1/4%
}{mr-1/2}}} \Biggr).
\]

\item $\operatorname{Gamma}(r,\lambda)$ (with $r$ known): $a=0$, $b=-\frac{1}{2r}$
and $%
H_{m}(X)=\sqrt{r}\ln({\frac{X}{rm-1/2}})$.

\item NEF--GHS$(r,\lambda)$ (with $r$ known): $a=0$, $b=-\frac
{1}{2r}$ and
\[
H_{m}(X)=\sqrt{r}\ln\Biggl( {\frac{X}{rm-1/2}}+\sqrt{1+{\frac{X^{2}}{%
(mr-1/2)^{2}}}} \Biggr).
\]
\end{itemize}

\begin{figure}[b]

\includegraphics[scale=0.99]{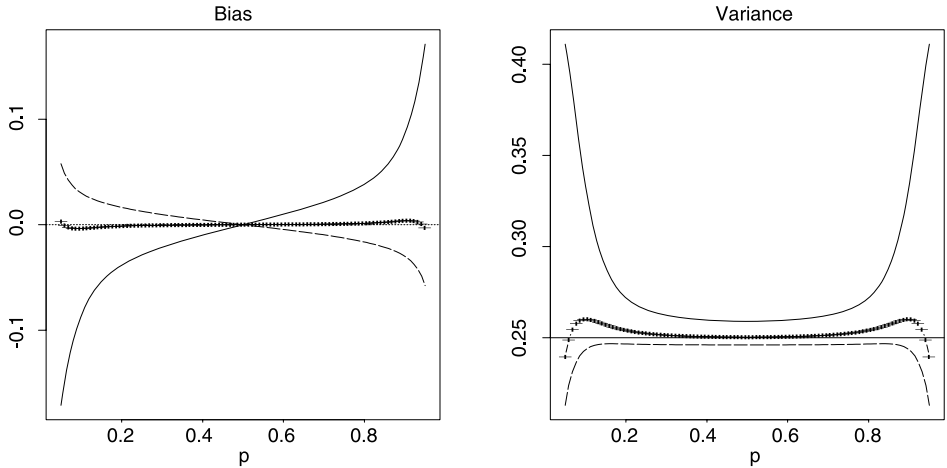}

\caption{Comparison of the mean (left panel) and variance
(right panel) of the arcsine transformations for $\operatorname{Binomial}(30, p)$ with
$c = 0$ (solid line), $c={\frac{1}{4}}$ ($+$ line) and $c={\frac{3}{8}}$ (dashed
line).}
\label{binomial-mean-var.plot}
\end{figure}

Note that the mean-matching VST is different from the more conventional VST
that optimally stabilizes the variance. Take the binomial distribution
with $%
r=1$ as an example. In this case the VST is an arcsine transformation.
Let $%
X_{1},\ldots,X_{m}\stackrel{\mathrm{i.i.d.}}{\sim}\operatorname{Bernoulli}(p)$ and then $%
X=\sum_{i=1}^{m}X_{i}\sim\operatorname{Binomial}(m,p)$. Figure \ref%
{binomial-mean-var.plot} compares the mean and variance of three arcsine
transformations of the form
\[
\arcsin\Biggl( \sqrt{\frac{X+c}{m+2c}} \Biggr)
\]
for the binomial variable $X$ with $m=30$. The choice of $c=0$ gives the
usual arcsine transformation, $c=3/8$ optimally stabilizes the variance
asymptotically, and $c=1/4$ yields the mean-matching arcsine transformation.
The left panel of Figure \ref{binomial-mean-var.plot} plots the bias
\[
\sqrt{m}\bigl(\mathbb{E}_{p}\arcsin\bigl(\sqrt{(X+c)/(m+2c)}\bigr)-\arcsin\bigl(\sqrt{p}\bigr)\bigr)
\]
as a function of $p$ for $c=0$, $c={\frac{1}{4}}$ and $c={\frac
{3}{8}}$. It
is clear from the plot that $c={\frac{1}{4}}$ is the best choice among the
three for matching the mean. On the other hand, the arcsine transformation
with $c=0$ yields significant bias and the transformation with
$c={\frac{3}{8%
}}$ also produces noticeably larger bias. The right panel plots the variance
of $\sqrt{m}\arcsin(\sqrt{(X+c)/(m+2c)})$ for $c=0$, $c={\frac
{1}{4}}$ and $%
c={\frac{3}{8}}$. Interestingly, over a wide range of values of $p$
near the
center the arcsine transformation with $c={\frac{1}{4}}$ is even slightly
better than the case with $c={\frac{3}{8}}$ and clearly $c=0$ is the worst
choice of the three. Figure \ref{poisson-mean-var.plot} below shows similar
behavior for the Poisson case.

\begin{figure}[b]

\includegraphics{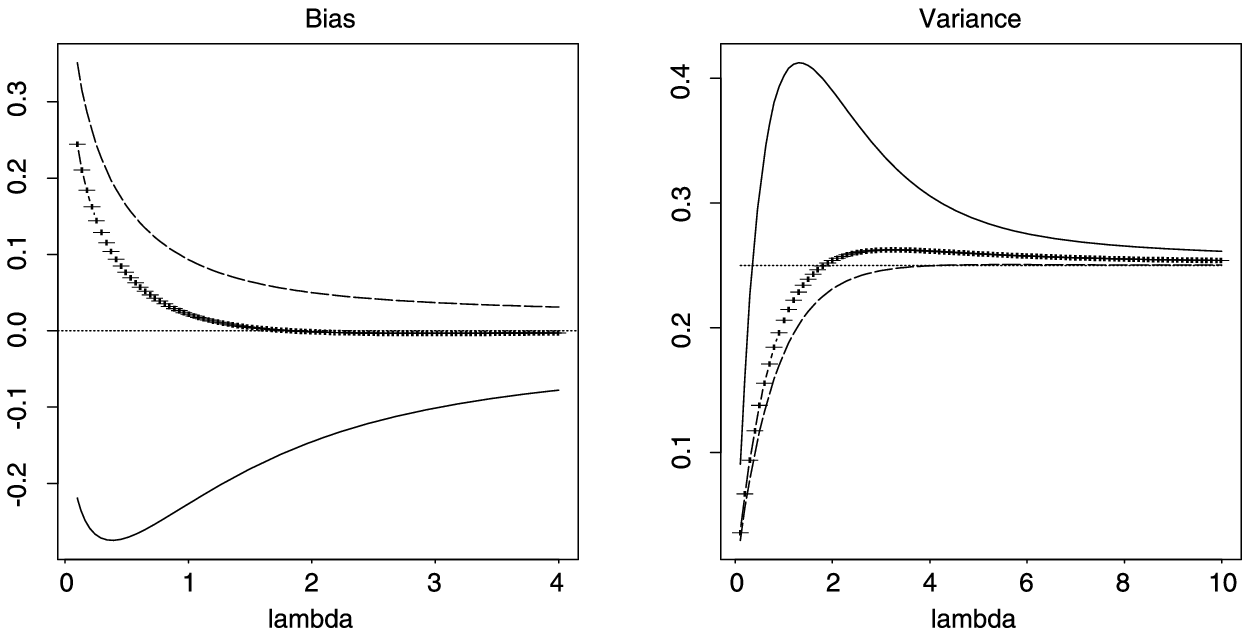}

\caption{Comparison\vspace*{1pt} of the mean (left panel) and variance
(right panel) of the root transformations for $\operatorname{Poisson}(\protect\lambda)$
with $c = 0$ (solid line), $c={\frac{1}{4}}$ ($+$ line) and $c={\frac{3}{8}}$
(dashed line).}
\label{poisson-mean-var.plot}
\end{figure}

\begin{figure}

\includegraphics{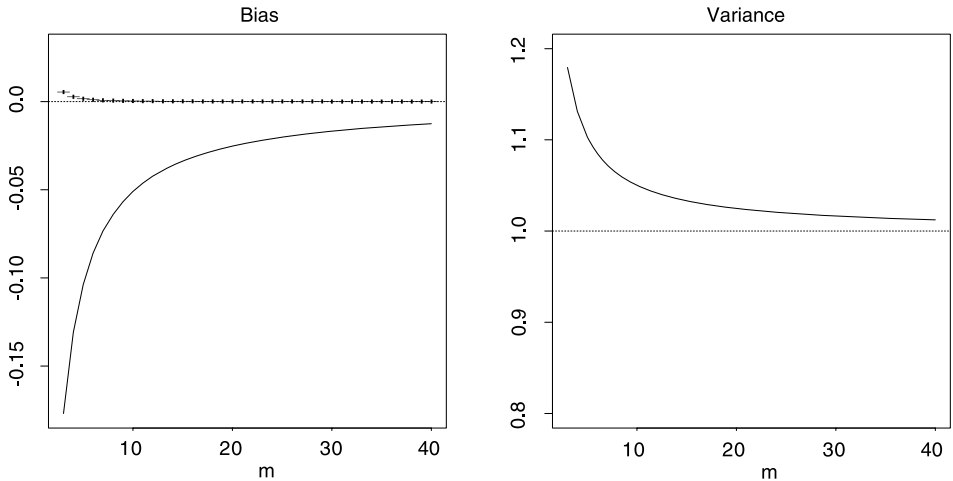}

\caption{Comparison of the mean (left panel) and variance
(right panel) of the log transformations for $\operatorname{Gamma}(m,\protect\lambda)$
with $c=0$ (solid line) and $c={\frac{1}{2}}$ ($+$ line).}
\label{gamma-mean-var.plot}
\end{figure}

Let us now consider the Gamma distribution with $r=1$ as an example for the
continuous case. The VST in this case is a log transformation. Let $%
X_{1},\ldots,X_{m}\stackrel{\mathrm{i.i.d.}}{\sim}\operatorname{Exponential}(\lambda)$.
Then $%
X=\sum_{i=1}^{m}X_{i}\sim\operatorname{Gamma}(m,\lambda)$. Figure \ref%
{gamma-mean-var.plot} compares the mean and variance of two log
transformations of the form
%
\begin{equation} \label{log.trans}
\ln\biggl( {\frac{X}{m-c}} \biggr)
\end{equation}
for the Gamma variable $X$ with $\lambda=1$ and $m$ ranging from 3 to 40.
The choice of $c=0$ gives the usual log transformation, and $c=1/2$ yields
the mean-matching log transformation. The left panel of Figure \ref%
{gamma-mean-var.plot} plots the bias as a function of $m$ for $c=0$ and
$c={%
\frac{1}{2}}$. It is clear from the plot that $c={\frac{1}{2}}$ is a much
better choice than $c=0$ for matching the mean. It is interesting to note
that in this case there do not exist constants $a$ and $b$ that optimally
stabilize the variance. The right panel plots the variance of $\sqrt
{m}\ln
(X)$, that is, $c=0$, as a function of $m$. In this case, it is obvious that
the variances are the same with $c=0$ and $c=1/2$ for the variable in~(\ref{log.trans}).
\begin{remark}
\label{remark1}
Mean-matching variance stabilizing transformations exist for some
other important families of distributions. We mention two that are
commonly used to model ``over-dispersed'' data. The first family is
often referred to as the gamma-Poisson family. See, for example,
Johnson, Kemp and Kotz (\citeyear{JKK05}), Berk and MacDonald (\citeyear{BM08}) and Hilbe (\citeyear{Hilbe07}).
Let $X_i|Zi \stackrel{\mathrm{ind}}{\sim} \operatorname{Poisson}(Z_i)$ with $Z_i
\stackrel{\mathrm{ind}}{\sim} \operatorname{Gamma}(\alpha, \sigma)$, $i = 1,\ldots
,m$. The
$Z_i$ are latent variables; only the $X_i$ are observed. The scale
parameter, $\sigma$, is assumed known, and the mean $\mu=\alpha
\sigma$
is the unknown parameter, $0 < \mu< \infty$. The resulting family of
distributions of each $X_i$ is a subfamily of the negative Binomial
$(r,p)$ with $p = (1+\sigma)^{-1}$, a fixed constant, and $r =
\mu/\sigma$. [Here this negative Binomial family is defined for
all $r > 0$ as having probability function,
$P(k ) = \Gamma(k + r) p^r (1-p)^k/\Gamma(k +1)\Gamma(r)$,
$k = 0,1,\ldots.$] This is a one-parameter family, but it is not an
exponential family. It can be verified that a mean-matching
variance stabilizing transformation for this family is given by
\[
Y=H_m(X)=2\sqrt{{X\over m} + {\sigma+ 1\over4m}}.
\]
This transformation has the desired properties (\ref{var}) and
(\ref{sharpbias}) with $G(\mu) = 2\sqrt{\mu}$. For the second family,
consider the beta-binomial family. See Johnson, Kemp and Kotz
(\citeyear{JKK05}). Here, $X_i| Z_i \stackrel{\mathrm{ind}}{\sim} \operatorname{Binomial} (r,Z_i )$ and
$Z_i \stackrel{\mathrm{ind}}{\sim} \operatorname{Beta} (a,b)$, $i =1,\ldots,m$.
Again, the $Z_i$ are latent variables; only the $X_i$ are observed.
For the family of interest here, we assume $a, b$ are allowed to vary
so that $a + b = k$, a known constant, and $0 < \mu= a/(a + b) <1$.
This family can alternatively be parameterized via
$\mu, \sigma= \mu(1-\mu)/(k+1)$. The resulting one-parameter
family of distributions of each $X_i$ is again not a one-parameter
exponential family. It can be verified that a mean-matching variance
stabilizing transformation for this family is given by
\[
Y=H_m(X)=2\sqrt{r} \arcsin\sqrt{X+ (\sigma+ 1)/4\over rm + (\sigma+1)/2}.
\]
This transformation has the desired properties (\ref{var}) and
(\ref{sharpbias}) with $G(\mu) = 2\times\break\arcsin\sqrt{\mu}$.
\end{remark}

\section{Nonparametric regression in exponential families}
\label{procedure.sec}

We now turn to nonparametric regression in exponential families. We begin
with the NEF--QVF. Suppose we observe
%
\begin{equation} \label{data1}
Y_i \stackrel{\mathrm{ind}}{\sim} \operatorname{NQ}(f(t_i)),\qquad i = 1,\ldots, n, t_i = {\frac
{i}{%
n}},
\end{equation}
and wish to estimate the mean function $f(t)$. In this setting, for the five
NEF--QVF families discussed in the last section the noise is not
additive and
non-Gaussian. Applying standard nonparametric regression methods
directly to
the data $\{Y_i\}$ in general do not yield desirable results. Our strategy
is to use the mean-matching VST to reduce this problem to a standard
Gaussian regression problem based on a sample $\{\widetilde Y_{j}\dvtx j = 1,
\ldots, T\}$ where
\[
\widetilde Y_{j} \sim N ( G ( f ( t_{j} ) )
,m^{-1} ) ,\qquad t_{j}=j/T, j=1,2,\ldots,T.
\]
Here $G$ is the VST defined in (\ref{vst}), $T$ is the number of bins,
and $%
m $ is the number of observations in each bin. The values of $T$ and $m$
will be specified later.

We begin by dividing the interval into $T$ equi-length subintervals
with $%
m=n/T$ observations in each subintervals. Let $Q_{j}$ be the sum of
observations on the $j$th subinterval $I_{j}=[{\frac{j-1}{T}},{\frac
{j}{T}}%
) $, $j=1,2,\ldots, T$,
%
\begin{equation} \label{Qj}
Q_{j}=\sum_{i=(j-1)m+1}^{jm}Y_{i}.
\end{equation}
The sums $\{Q_{j}\}$ can be treated as observations for a Gaussian
regression directly, but this in general leads to a heteroscedastic problem.
Instead, we apply the mean-matching VST discussed in Section \ref{vst.sec},
and then treat $H_{m}(Q_{j}) $ as new observations in a homoscedastic
Gaussian regression problem. To be more specific, let
%
\begin{equation}\label{binned.trans.data}
Y_{j}^{\ast}=H_m(Q_j)=G\biggl(\frac{Q_{j}+a}{m+b}\biggr),\qquad j=1,\ldots,T,
\end{equation}
where the constants ${a}$ and $b$ are chosen as in (\ref%
{optimalconstants}) to match the means. The transformed data $Y^{\ast
}=(Y_{1}^{\ast},\ldots,Y_{T}^{\ast})$ is then treated as the new
equi-spaced sample for a nonparametric Gaussian regression problem.

We will first estimate $G(f(t_{i}))$, then take a transformation of the
estimator to estimate the mean function $f$. After the original regression
problem is turned into a Gaussian regression problem through binning
and the
mean-matching VST, in principle any good nonparametric Gaussian regression
method can be applied to the transformed data $\{Y_{j}^*\}$ to
construct an
estimate of $G(f(\cdot))$. The general ideas for our approach can be
summarized as follows.

\begin{enumerate}
\item \textit{Binning}: divide $\{Y_i\}$ into $T$ equal length
intervals between 0 and 1. Let $Q_1, Q_2,\ldots, Q_T$ be the sum of the
observations in each of the intervals. Later results suggest a choice
of $T$
satisfying $T\asymp n^{3/4}$ for the NEF--QVF case and $T\asymp
n^{1/2}$ for
the non-QVF case. See Section \ref{theory.sec} for details.

\item \textit{VST}: let $Y_{j}^*=H_m(Q_j)$, $j=1,\ldots
,T$, and
treat $Y^*=(Y_{1}^*,Y_{2}^*,\ldots,Y_{T}^*)$ as the new equi-spaced sample
for a nonparametric Gaussian regression problem.

\item \textit{Gaussian regression}: apply your favorite
nonparametric regression procedure to the binned and transformed data $Y^*$
to obtain an estimate $\widehat{G ( f ) }$ of $G ( f ) $.

\item \textit{Inverse VST}: estimate the mean function
$f$ by $%
\widehat{f}=G^{-1} ( \widehat{G ( f ) } ) $. If $\widehat{%
G ( f ) }$ is not in the domain of $G^{-1}$ which is an interval
between $a$ and $b$ ($a$ and $b$ can be $\infty$), we set $G^{-1} (
\widehat{G(f) } ) =G^{-1}(a) $ if $\widehat{G(f)} < a$ and set $%
G^{-1} ( \widehat{G( f) } ) = G^{-1} ( b ) $ if $\widehat{%
G ( f ) } > b$. For example, $G^{-1} ( a ) =0$ when $a<0$
in the case of negative Binomial and NEF--GHS distributions.
\end{enumerate}

\subsection{Effects of binning and VST}
\label{effect.sec}

As mentioned earlier, after binning and the mean-matching VST, one can treat
the transformed data $\{Y_j^*\}$ as if they were data from a homoscedastic
Gaussian nonparametric regression problem. A key step in understanding why
this procedure works is to understand the effects of binning and the VST.
Quantile coupling provides an important technical tool to shed insights on
the procedure.

The following result, which is a direct consequence of the quantile coupling
inequality of Koml\'{o}s, Major and Tusn\'{a}dy (\citeyear{KMT75}), shows that the
binned and transformed data can be well approximated by independent normal
variables.
\begin{lemma}
\label{KMT.lem} Let $X_{i}\stackrel{\mathrm{i.i.d.}}{\sim}\operatorname{NQ}(\mu)$ with
variance $V$
for $i=1,\ldots,m$ and let $X=\sum_{i=1}^{m}X_{i}$. Under the
assumptions of
Lemma \ref{approx.lem}, there exists a standard normal random variable
$%
Z\sim N(0,1)$ and constants $c_{1},c_{2},c_{3}>0$ not depending on $m$ such
that whenever the event $A=\{|X-m\mu|\leq c_{1}m\}$ occurs,
%
\begin{equation} \label{KMT}
\bigl|X-m\mu-\sqrt{mV}Z\bigr|<c_{2}Z^{2}+c_{3}.
\end{equation}
\end{lemma}

Hence, for large $m$, $X$ can be treated as a normal random variable with
mean $m\mu$ and variance $mV$. Let $Y=H_{m}(X)=G({\frac{X+a}{m+b}})$,
$%
\epsilon=\mathbb{E}Y-G({\mu})$ and $Z$ be a standard normal variable
satisfying (\ref{KMT}). Then $Y$ can be written as
%
\begin{equation}\label{Y.decomp0}
Y=G({\mu})+\epsilon+m^{-{{1/2}}}Z+\xi,
\end{equation}
where
%
\begin{equation}\label{xi}
\xi=G\biggl({\frac{X+a}{m+b}}\biggr)-G({\mu})-\epsilon-m^{-{{1/2}}}Z.
\end{equation}
In the decomposition (\ref{Y.decomp0}), $\epsilon$ is the deterministic
approximation error between the mean of $Y$ and its target value
$G({\mu})$
and $\xi$ is the stochastic error measuring the difference of $Y$ and its
normal approximation. It follows from Lemma \ref{approx.lem} that when $m$
is large, $\epsilon$ is ``small,''
$|\epsilon|\leq{c}m^{-2}$ for some constant $c>0$. The following result,
which is proved in Section \ref{coupling-proof.sec}, shows that the random
variable $\xi$ is ``stochastically
small.''
\begin{lemma}
\label{xi.lem} Let $X_{i}\stackrel{\mathrm{i.i.d.}}{\sim}\operatorname{NQ}(\mu)$ with variance
$V$ for
$i=1,\ldots,m$, and $X=\sum_{i=1}^{m}X_{i}$. Let $Z$ be the standard normal
variable given as in Lemma \ref{KMT.lem} and let $\xi$ be given as in
(\ref%
{xi}). Then for any integer $k\geq1$ there exists a constant $C_{k}>0$ such
that for all $\lambda\geq1$ and all $a>0$,
%
\begin{equation} \label{moment.bd}
\mathbb{E}|\xi|^{k}\leq C_{k}m^{-k} \quad\mbox{and}\quad \mathbb{P}(|\xi
|>a)\leq C_{k}(am)^{-k}.
\end{equation}
\end{lemma}

The discussion so far has focused on the effects of the VST for i.i.d.
observations. In the nonparametric function estimation problem mentioned
earlier, observations in each bin are independent but not identically
distributed since the mean function $f$ is not a constant in general.
However, through coupling, observations in each bin can in fact be treated
as if they were i.i.d. random variables when the function $f$ is
smooth. Let
$X_{i}\sim \operatorname{NQ}(\mu_{i})$, $i=1,\ldots,m$, be independent. Here the means
$\mu
_{i}$ are ``close'' but not equal. Let $\mu
$ be a value close to the $\mu_{i}$'s. The analysis given in Section
\ref{coupling-proof.sec} shows that $X_{i}$ can in fact be coupled with i.i.d.
random variables $X_{i,c}$ where $X_{i,c}\stackrel{\mathrm{i.i.d.}}{\sim}\operatorname{NQ}(\mu
)$. See
Lemma \ref{iidapprox} in Section \ref{coupling-proof.sec} for a precise
statement.

How well the transformed data $\{Y_{j}^{\ast}\}$ can be approximated
by an
ideal Gaussian regression model depends partly on the smoothness of the mean
function $f$. For $0<d\leq1$, define the Lipschitz class $\Lambda^{d}(M)$
by
\[
\Lambda^{d}(M) = \{f\dvtx|f(t_{1})-f(t_{2})|\leq
M |t_{1}-t_{2}|^{d} 0\leq t_{1}, t_{2}\leq1\}
\]
and%
\[
F^{d}(M,\varepsilon,v)=\{f\dvtx f\in\Lambda^{d}(M),f(t)\in[ \varepsilon
,v ] , \mbox{for all
$x\in[0,1]$}\},
\]
where $[\varepsilon,v]$ with $\epsilon<v$ is a compact set in the interior
of the mean parameter space of the natural exponential family. Lemmas
\ref%
{approx.lem}, \ref{KMT.lem}, \ref{xi.lem} and \ref{iidapprox}
together yield
the following result which shows how far away are the transformed data $
\{Y_{j}^{\ast}\}$ from the ideal Gaussian model.
\begin{theorem}
\label{Yi.thm} Let $Y_{j}^{\ast}=G(\frac{Q_{j}+a}{m+b})$ be given as
in $(%
\ref{binned.trans.data})$ and let $f\in F^{d}(M,\varepsilon,v)$. Then
$%
Y_{j}^{\ast}$ can be written as
%
\begin{equation} \label{Yi.decomp}
Y_{j}^{\ast}=G\biggl(f\biggl({\frac{j}{T}}\biggr)\biggr)+\epsilon_{j}+m^{-{
{1/2}}}Z_{j}+\xi
_{j},\qquad j=1,2,\ldots,T,
\end{equation}
where $Z_{j}\stackrel{\mathrm{i.i.d.}}{\sim}N(0,1)$, $\epsilon_{j}$ are constants
satisfying $|\epsilon_{j}|\leq{c} ( m^{-2}+T^{-d} ) $ and
consequently for some constant $C>0$
%
\begin{equation}\label{epsilon.bnd}
{\frac{1}{T}}\sum_{j=1}^{T}\epsilon_{j}^{2}\leq C (
m^{-4}+T^{-2d} )
\end{equation}
and $\xi_{j}$ are independent and ``stochastically
small'' random variables satisfying that for any integer $%
k>0$ and any constant $a>0$
%
\begin{eqnarray} \label{moment.bd1}
\mathbb{E}|\xi_{j}|^{k} &\leq& C_{k}\log^{2k}m\cdot(m^{-k}+T^{-dk}) %
\quad\mbox{and}\nonumber\\[-8pt]\\[-8pt]
\mathbb{P}(|\xi_{j}|>a)&\leq& C_{k}\log^{2k}m\cdot
(m^{-k}+T^{-dk})a^{-k},\nonumber
\end{eqnarray}
where $C_{k}>0$ is a constant depending only on $k,d$ and $M$.
\end{theorem}

Theorem \ref{Yi.thm} provides explicit bounds for both the
deterministic and
stochastic errors. This is an important technical result which serves
as a
major tool for the proof of the main results given in Section \ref%
{theory.sec}.
\begin{remark}
There is a tradeoff between the two terms in the bound (\ref%
{epsilon.bnd}) for the overall approximation error ${\frac{1}{T}}%
\sum_{j=1}^{T}\epsilon_{j}^{2}$. There are two sources to the approximation
error: one is the variation of the functional values within a bin and one
comes from\vspace*{-1pt} the expansion of the mean of $Y_{j}^{\ast}$ (see Lemma \ref
{approx.lem}). The former is related to the smoothness of the function $f$
and is controlled by the $T^{-2d}$ term and the latter is bounded by
the $%
m^{-4}$ term. In addition, there is the discretization error between the
sampled function $\{G(f({j/T}))\dvtx j=1,\ldots,T\}$ and the whole function
$G(f(t))$%
, which is obviously a decreasing function of $T$. Furthermore, the choice
of $T$ also affects the stochastic error $\xi$. A good choice of $T$ makes
all three types of errors negligible relative to the minimax risk. See
Section \ref{discussion.sec} for further discussions.
\end{remark}
\begin{remark}
In Section \ref{theory.sec} we introduce Besov balls $%
B_{p,q}^{\alpha}(M)$ for the analysis of wavelet regression methods. A Besov
ball $B_{p,q}^{\alpha}(M)$ can be embedded into a Lipschitz class
$\Lambda
^{d} ( M') $ with $d=\min( \alpha-1/p,1 ) $ and some $M'>0$.
\end{remark}

Although the main focus of this paper is on the NEF--QEF, our method of
binning and VST can be extended to the general one-parameter NEF. This
extension is discussed in Section \ref{General-NEF.sec} where a
version of
Theorem \ref{Yi.thm} for the standard VST is developed in the general case.

\subsection{Wavelet thresholding}
\label{waveletblock.sec}

One can apply any good nonparametric Gaussian regression procedure to the
transformed data $\{Y_{j}^*\}$ to construct an estimator of the
function $f$%
. To illustrate our general methodology, in the present paper we shall use
wavelet block thresholding to construct the final estimators of the
regression function. Before we can give a detailed description of our
procedures, we need a brief review of basic notation and definitions.

Let $\{\phi, \psi\}$ be a pair of father and mother wavelets. The
functions $\phi$ and $\psi$ are assumed to be compactly supported and
$%
\int\phi=1$, and dilation and translation of $\phi$ and $\psi$ generates
an orthonormal wavelet basis. For simplicity in exposition, in the present
paper we work with periodized wavelet bases on $[0,1]$. Let
\[
\phi_{j,k}^{p}(t)=\sum_{l=-\infty}^{\infty}\phi_{j,k}(t-l),\qquad \psi
_{j,k}^{p}(t)=\sum_{l=-\infty}^{\infty}\psi_{j,k}(t-l)\qquad
\mbox{for $t \in[0,1]$},
\]
where $\phi_{j,k}(t)=2^{j/2}\phi(2^{j}t-k)$ and $\psi
_{j,k}(t)=2^{j/2}\psi(2^{j}t-k)$. The collection \{$\phi
_{j_{0},k}^{p}, k=1,\ldots,2^{j_{0}}; \psi_{j,k}^{p}, j\geq j_{0}\geq
0,k=1,\ldots,2^{j}$\} is then an orthonormal basis of $L^{2}[0,1]$, provided
the primary resolution level $j_{0}$ is large enough to ensure that the
support of the scaling functions and wavelets at level $j_{0}$ is not the
whole of $[0,1]$. The superscript ``$p$''
will be suppressed from the notation for convenience. An orthonormal wavelet
basis has an associated orthogonal Discrete Wavelet Transform (DWT) which
transforms sampled data into the wavelet coefficients. See Daubechies (\citeyear{Daubechies92})
and Strang (\citeyear{strang92}) for further details about the wavelets and discrete
wavelet transform. A square-integrable function $f$ on $[0,1]$ can be
expanded into a wavelet series:
%
\begin{equation} \label{resolution}
f(t)=\sum_{k=1}^{2^{j_{0}}}\widetilde{\theta}_{j_{0},k}\phi
_{j_{0},k}(t)+\sum_{j=j_{0}}^{\infty}\sum_{k=1}^{2^{j}}\theta
_{j,k}\psi
_{j,k}(t),
\end{equation}
where $\widetilde{\theta}_{j,k}=\langle f,\phi_{j,k}\rangle, \theta
_{j,k}=\langle f,\psi_{j,k}\rangle$ are the wavelet coefficients of $f$.

\subsection{Wavelet procedures for generalized regression}
\label{VST.procedure.sec}

We now give a detailed description of the wavelet thresholding procedures
BlockJS and NeighCoeff in this section and study the properties of the
resulting estimators in Section \ref{theory.sec}. We shall show that our
estimators enjoy a high degree of adaptivity and spatial adaptivity and are
easily implementable.

Apply the discrete wavelet transform to the binned and transformed data
$Y^*$%
, and let $U=T^{-{{1/2}}}WY^*$ be the empirical wavelet coefficients,
where $W$ is the discrete wavelet transformation matrix. Write
%
\begin{equation} \label{wavelet.coeff}\qquad
U=(\widetilde{y}_{j_{0},1},\ldots,\widetilde{y}%
_{j_{0},2^{j_{0}}},y_{j_{0},1},\ldots,y_{j_{0},2^{j_{0}}},\ldots
,y_{J-1,1},\ldots,y_{J-1,2^{J-1}})^{\prime}.
\end{equation}

Here $\widetilde{y}_{j_{0},k}$ are the gross structure terms at the lowest
resolution level, and $y_{j,k}$ ($j=j_{0},\ldots,J-1,k=1,\ldots,2^{j}$)
are empirical wavelet coefficients at level $j$ which represent fine
structure at scale $2^{j}$. The empirical wavelet coefficients can then be
written as
%
\begin{equation}\label{ujk.decomp}
y_{j,k}=\theta_{j,k}+\epsilon_{j,k}+{\frac{1}{\sqrt{n}}}z_{j,k}+\xi_{j,k},
\end{equation}
where $\theta_{j,k}$ are the true wavelet coefficients of $G(f)$,
$\epsilon
_{j,k}$ are ``small'' deterministic
approximation errors, $z_{j,k}$ are i.i.d. $N(0,1)$, and $\xi_{j,k}$ are
some ``small'' stochastic errors. The
theoretical calculations given in Section \ref{proof.sec} will show that
both $\epsilon_{j,k}$ and $\xi_{j,k}$ are negligible.
If these negligible errors are ignored then we have
%
\begin{equation}\label{idealized.decomp}
y_{j,k}\approx\theta_{j,k}+{\frac{1}{\sqrt{n}}}z_{j,k},
\end{equation}
which is the idealized Gaussian sequence model with noise level $\sigma
=1/%
\sqrt{n}$. Both BlockJS [Cai (\citeyear{Cai99})] and NeighCoeff [Cai and Silverman (\citeyear{CS01})]
were originally developed for this ideal model. Here we shall apply these
methods to the empirical coefficients $y_{j,k}$ as if they were
observed as
in (\ref{idealized.decomp}).

\begin{figure}

\includegraphics{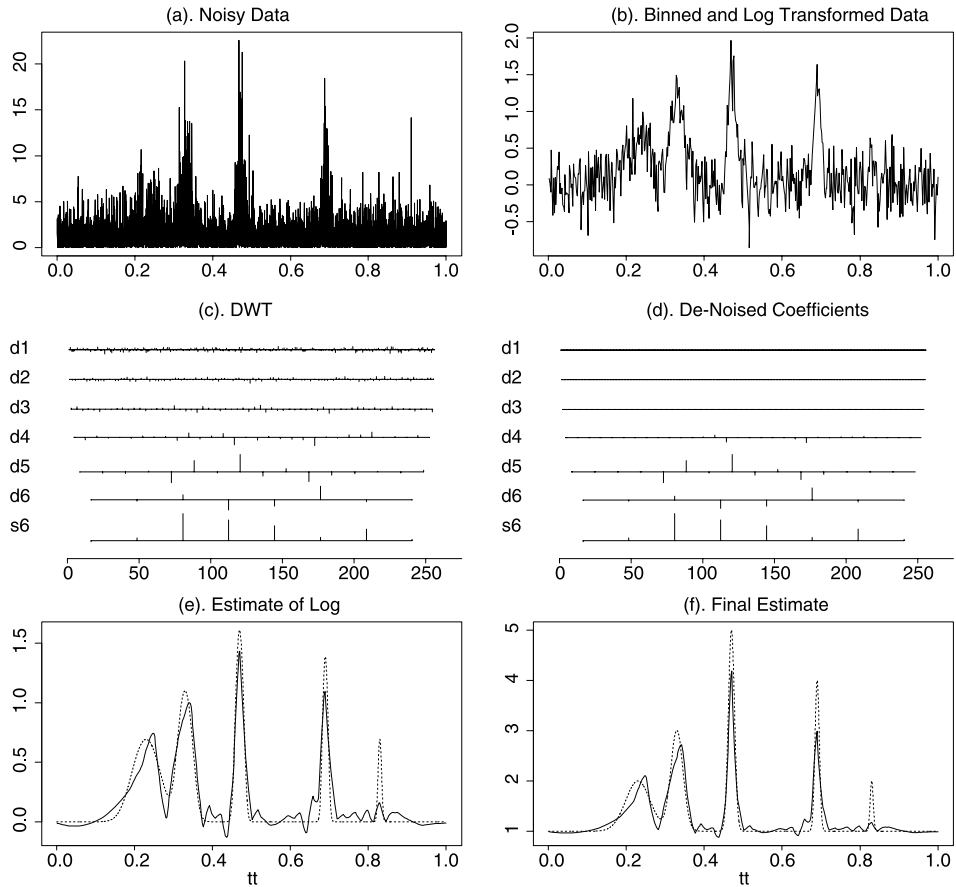}

\caption{An example of nonparametric Gamma regression using
the mean-matching VST and wavelet block thresholding.}
\label{Spike.Gamma.ex}
\end{figure}

We first describe the \textit{BlockJS} procedure. At each resolution
level $%
j $, the empirical wavelet coefficients $y_{j,k}$ are grouped into
nonoverlapping blocks of length $L$. As in the sequence estimation setting
let $B_{j}^{i}=\{(j,k)\dvtx(i-1)L+1\leq k\leq iL\}$ and let
$S_{j,i}^{2}\equiv
\sum_{(j,k)\in B_{j}^{i}}y_{j,k}^{2}$. A modified James--Stein
shrinkage rule
is then applied to each block $B_{j}^{i}$, that is,
%
\begin{equation} \label{block.est}
\widehat{\theta}_{j,k}= \biggl( 1-{\frac{\lambda_{\ast}L}{nS_{j,i}^{2}}} \biggr)
_{+} y_{j,k} \qquad\mbox{for $(j,k)\in
B_{j}^{i}$},
\end{equation}
where $\lambda_{\ast}=4.50524$ is the solution to the equation
$\lambda
_{\ast}-\log\lambda_{\ast}=3$ [see Cai (\citeyear{Cai99}) for details], and
$\frac{1%
}{n}$ is approximately the variance of each $y_{j,k}$. For the gross
structure terms at the lowest resolution level $j_{0}$, we set $\widehat{\!\widetilde{
\theta}}_{j_{0},k}=\widetilde{y}_{j_{0},k}$. The estimate of $G(f(\cdot
))$ at
the equally spaced sample points $\{{\frac{i}{T}}\dvtx i=1,\ldots,T\}$ is then
obtained by applying the inverse discrete wavelet transform (IDWT) to the
denoised wavelet coefficients. That is, $\{G(f({\frac
{i}{T}}))\dvtx i=1,\ldots
,T\}$ is estimated by $\widehat{G(f)}=\{\widehat{G(f({\frac{i}{T}}))}
\dvtx i=1,\ldots,T\}$ with $\widehat{G(f)}=T^{{1/2}}W^{-1}\cdot\widehat
{\theta%
}$. The estimate of the whole function $G(f)$ is given by
\[
\widehat{G(f(t))}=\sum_{k=1}^{2^{j_{0}}}\widehat{\!\widetilde{\theta
}}_{j_{0},k}\phi
_{j_{0},k}(t)+\sum_{j=j_{0}}^{J-1}\sum_{k=1}^{2^{j}}\widehat{\theta
}_{j,k}\psi
_{j,k}(t).
\]
The mean function $f$ is estimated by
%
\begin{equation} \label{fun.est}
\widehat{f}_{\mathrm{BJS}}(t)=G^{-1}(\widehat{G(f(t))}).
\end{equation}
Figure \ref{Spike.Gamma.ex} shows the steps of the procedure for an example
in the case of nonparametric Gamma regression.

We now turn to the \textit{NeighCoeff} procedure. This procedure, introduced
in Cai and Silverman (\citeyear{CS01}) for Gaussian regression, incorporates
information about neighboring coefficients in a different way from the
BlockJS procedure. NeighCoeff also thresholds the empirical
coefficients in
blocks, but estimates wavelet coefficients individually. It chooses a
threshold for each coefficient by referencing not only to that coefficient
but also to its neighbors. As shown in Cai and Silverman (\citeyear{CS01}), NeighCoeff
outperforms BlockJS numerically, but with slightly inferior asymptotic
properties.

Let the empirical coefficients $\{y_{j,k}\}$ be given same as before. To
estimate a coefficient $\theta_{j,k}$ at resolution level $j$, we form a
block of size $3$ by including the coefficient $y_{j,k}$ together with its
immediate neighbors $y_{j,k-1}$ and $y_{j,k+1}$. (If periodic boundary
conditions are not being used, then for the two coefficients at the boundary
blocks, again of length $3$, are formed by only extending in one direction.)
Estimate the coefficient $\theta_{j,k}$ by
%
\begin{equation}
\widehat{\theta}_{j,k}= \biggl( 1-{\frac{2\log n}{nS_{j,k}^{2}}} \biggr)
_{+}y_{j,k},
\end{equation}
where $S_{j,k}^{2}=y_{j,k-1}^{2}+y_{j,k}^{2}+y_{j,k+1}^{2}$. The gross
structure terms at the lowest resolution level are again estimated by
$\widehat{\!\widetilde{\theta}}_{j_{0},k}=\widetilde{y}_{j_{0},k}$. The rest of the
steps are
same as before. Namely, the inverse DWT is applied to obtain an
estimate $%
\widehat{G(f)}$ and the mean function $f$ is then estimated by
$\widehat{f}%
_{\mathrm{NC}}(t)=G^{-1}(\widehat{G(f(t))})$.

We can envision a sliding window of size $3$ which moves one position each
time and only the middle coefficient in the center is estimated for a given
window. Each individual coefficient is thus shrunk by an amount that depends
on the coefficient and on its immediate neighbors. Note that NeighCoeff uses
a lower threshold level than the universal thresholding procedure of Donoho
and Johnstone (\citeyear{DJ94}). In NeighCoeff, a coefficient is estimated by zero only
when the sum of squares of the empirical coefficient and its immediate
neighbors is less than $2\sigma^{2}\log n$, or the average of the squares
is less than ${\frac{2}{3}}\sigma^{2}\log n$.

\section{Theoretical properties}
\label{theory.sec}

In this section we investigate the asymptotic properties of the procedures
proposed in Section \ref{procedure.sec}. Numerical results will be
given in Section \ref{numerical.sec}.

We study the theoretical properties of our procedures over the Besov spaces
that are by now standard for the analysis of wavelet regression methods.
Besov spaces are a very rich class of function spaces and contain as special
cases many traditional smoothness spaces such as H\"{o}lder and Sobolev
spaces. Roughly speaking, the Besov space $B_{p,q}^{\alpha}$ contains
functions having $\alpha$ bounded derivatives in $L^{p}$ norm, the third
parameter $q$ gives a finer gradation of smoothness. Full details of Besov
spaces are given, for example, in Triebel (\citeyear{triebel92}) and DeVore and Popov
(\citeyear{DP88}). A wavelet $\psi$ is called \textit{r-regular} if $\psi$ has $r$
vanishing moments and $r$ continuous derivatives. For a given r-regular
mother wavelet $\psi$ with $r>\alpha$ and a fixed primary resolution level
$j_{0}$, the Besov sequence norm $\Vert\cdot\Vert_{b_{p,q}^{\alpha
}}$ of
the wavelet coefficients of a function $f$ is then defined by
%
\begin{equation} \label{besseqnorm}
\Vert f\Vert_{b_{p,q}^{\alpha}}=\Vert\underline{\xi}_{j_{0}}\Vert
_{p}+ \Biggl( \sum_{j=j_{0}}^{\infty}(2^{js}\Vert\underline{\theta}%
_{j}\Vert_{p})^{q} \Biggr) ^{{1/q}},
\end{equation}
where $\underline{\xi}_{j_{0}}$ is the vector of the father wavelet
coefficients at the primary resolution level $j_{0}$, $\underline
{\theta}%
_{j}$ is the vector of the wavelet coefficients at level $j$, and
$s=\alpha+%
{\frac{1}{2}}-{\frac{1}{p}}>0$. Note that the Besov function norm of
index $%
(\alpha,p,q)$ of a function $f$ is equivalent to the sequence norm
(\ref%
{besseqnorm}) of the wavelet coefficients of the function. See Meyer (\citeyear{meyer92}).
We define%
%
\begin{equation}\label{Besovball}
B_{p,q}^{\alpha} ( M ) = \{ f;\Vert f\Vert_{b_{p,q}^{\alpha
}}\leq M \}
\end{equation}
and
%
\begin{equation}
F_{p,q}^{\alpha}(M,\varepsilon,v)=\{f\dvtx f\in B_{p,q}^{\alpha
}(M),f(t)\in
[\varepsilon,v] \mbox{ for all $t\in[0,1]$}\},
\end{equation}
where $[\varepsilon,v]$ with $\epsilon<v$ is a compact set in the interior
of the mean parameter space of the natural exponential family.

The following theorem shows that our estimators achieve near optimal global
adaptation under integrated squared error for a wide range of Besov balls.
\begin{theorem}
\label{global.thm} Suppose the wavelet $\psi$ is r-regular. Let
$X_{i}\sim
\operatorname{NQ}(f(t_{i})), i=1,\ldots,n, t_{i}={\frac{i}{n}}$. Let $T=cn^{{3/4}}$.
Then the estimator $\widehat{f}_{\mathrm{BJS}}$ defined in (\ref{fun.est}) satisfies
\begin{eqnarray*}
\sup_{f\in F_{p,q}^{\alpha}(M,\varepsilon,v)}\mathbb{E}\Vert\widehat
{f}%
_{\mathrm{BJS}}-f\Vert_{2}^{2}
\leq\cases{
Cn^{-{({2\alpha})/({1+2\alpha})}}, \cr\qquad \mbox{$p\geq2, \alpha\leq r$
and ${\dfrac{3}{2}\biggl(\alpha-{\dfrac{1}{p}}\biggr)>{\dfrac{2\alpha}{1+2\alpha}}}$},
\vspace*{2pt}\cr
Cn^{-{({2\alpha})/({1+2\alpha})}}(\log n)^{({2-p})/({p(1+2\alpha
)})},\cr \qquad\mbox{$1\leq p<2, \alpha\leq r$ and ${\dfrac{3}{2}\biggl(\alpha-{\dfrac
{1}{p}}\biggr)>{%
\dfrac{2\alpha}{1+2\alpha}}}$},}
\end{eqnarray*}
and the estimator $\widehat{f}_{\mathrm{NC}}$ satisfies
\begin{eqnarray}
\sup_{f\in F_{p,q}^{\alpha}(M,\varepsilon,v)}\mathbb{E}\Vert\widehat
{f}%
_{\mathrm{NC}}-f\Vert_{2}^{2}\leq C \biggl( \frac{\log n}{n} \biggr) ^{{({2\alpha})/({%
1+2\alpha})}},\nonumber\\
\eqntext{p\geq1, \alpha\leq r\mbox{ and }{\dfrac{3}{2}\biggl(\alpha
-{%
\dfrac{1}{p}}\biggr)>{\dfrac{2\alpha}{1+2\alpha}}}.}
\end{eqnarray}
\end{theorem}
\begin{remark}
Note that when $f(t)\in[\varepsilon,v]$, the condition
$f\in
B_{p,q}^{\alpha}(M)$ implies that there exists $M^{\prime}>0$ such
that $%
G ( f ) \in B_{p,q}^{\alpha} ( M^{\prime} ) $ with
\[
M^{\prime}=c_{0}+cM \Biggl[ \sum_{l=1}^{ \lfloor\alpha\rfloor
+1}c_{l}v^{l-1}+c_{ \lfloor\alpha\rfloor+1} \Biggr] \qquad
\mbox{for some }c>0,
\]
where $c_{l}=\sup_{y\in[\varepsilon,v]} \vert G^{ (
l ) } ( y ) \vert$ with $l=0,\ldots, \lfloor
\alpha\rfloor+1$, since it follows from Theorem 3 on page 344 and
Remark 3
on page 345 of Runst (\citeyear{Runst86}) that
\begin{eqnarray*}
\Vert G ( f ) \Vert_{B_{p,q}^{\alpha}} &\leq& \Vert G ( f )
\Vert_{p}\\
&&{} + c\Vert f\Vert_{B_{p,q}^{\alpha}} \Biggl[ \sum_{l=1}^{ \lfloor
\alpha\rfloor+1} \bigl\Vert G^{ ( l ) } ( f )
\bigr\Vert_{\infty} \Vert f \Vert_{\infty}^{l-1}+ \bigl\Vert
G^{ \lfloor\alpha\rfloor+1} ( f ) \bigr\Vert_{\infty
} \Biggr].
\end{eqnarray*}
\end{remark}
\begin{remark}
Simple algebra\vspace*{-2pt} shows that $\frac{3}{2}(\alpha-{\frac
{1}{p}})>{\frac{%
2\alpha}{1+2\alpha}}$ is equivalent to ${\frac{2\alpha^{2}-\alpha
/3}{%
1+2\alpha}}>{\frac{1}{p}}$. This\vspace*{1pt} condition is needed to ensure that the
discretization error over the Besov ball $B_{p,q}^{\alpha}(M)$ is
negligible\vspace*{1pt} relative to the minimax risk. See Section \ref{discussion.sec}
for more discussions.
\end{remark}

For functions of spatial inhomogeneity, the local smoothness of the
functions varies significantly from point to point and global risk
given in
Theorem \ref{global.thm} cannot wholly reflect the performance of estimators
at a point. We use the local risk measure
%
\begin{equation}\label{l.risk}
R(\widehat{f}(t_{0}), f(t_{0}))=\mathbb{E}\bigl(\widehat{f}(t_{0})-f(t_{0})\bigr)^{2}
\end{equation}
for spatial adaptivity.

The local smoothness of a function can be measured by its local H\"{o}lder
smoothness index. For a fixed point $t_{0}\in(0,1)$ and $0<\alpha\leq1$,
define the local H\"{o}lder class $\Lambda^{\alpha}(M,t_{0},\delta)$ as
follows:
\[
\Lambda^{\alpha}(M,t_{0},\delta) = \{f\dvtx|f(t)-f(t_{0})|\leq
M |t-t_{0}|^{\alpha}, \mbox{for $t\in(t_0-\delta, t_0+\delta)$}\}.
\]
If $\alpha>1$, then
\[
\Lambda^{\alpha}(M,t_{0},\delta) = \bigl\{f\dvtx\bigl|f^{(\lfloor\alpha\rfloor
)}(t)-f^{(\lfloor\alpha\rfloor)}(t_{0})\bigr|\leq M |t-t_{0}|^{\alpha
^{\prime}} \mbox{ for $t\in(t_0-\delta,
t_0+\delta)$}\bigr\},
\]
where $\lfloor\alpha\rfloor$ is the largest integer less than
$\alpha$
and $\alpha^{\prime}=\alpha-\lfloor\alpha\rfloor$. Define%
\[
F^{\alpha}(M,t_{0},\delta,\varepsilon,v)=\{f\dvtx f\in\Lambda^{\alpha
}(M,t_{0},\delta),f(x)\in[ \varepsilon,v ] %
\mbox{ for all
$x\in[0,1]$}\}.
\]

In Gaussian nonparametric regression setting, it is a well-known fact that
for estimation at a point, one must pay a price for adaptation. The optimal
rate of convergence for estimating $f(t_0)$ over function class $%
\Lambda^\alpha(M,t_0,\delta)$ with $\alpha$ completely known is $n^{-2
\alpha/(1 + 2 \alpha)}$. Lepski (\citeyear{Lepski90}) and Brown and Low (\citeyear{BL96})
showed that
one has to pay a price for adaptation of at least a logarithmic factor. It
is shown that the local adaptive minimax rate over the H\"{o}lder class
$%
\Lambda^\alpha(M, t_0, \delta)$ is $(\log n/n)^{2 \alpha/(1 + 2
\alpha)}$.

The following theorem shows that our estimators achieve optimal local
adaptation with the minimal cost.
\begin{theorem}
\label{LocalAdapt.thm} Suppose the wavelet $\psi$ is r-regular with $%
1/6<\alpha\leq r$. Let $t_{0}\in(0,1)$ be fixed. Let $X_{i}\sim
\operatorname{NQ}(f(t_{i})), i=1,\ldots,n, t_{i}=\frac{i}{n}$. Let $T=cn^{{3/4}}$.
Then for $\widehat{f}=\widehat{f}_{\mathrm{BJS}}$ or $\widehat{f}_{\mathrm{NC}}$
%
\begin{equation} \label{point}
\sup_{F^{\alpha}(M,t_{0},\delta,\varepsilon,v)}\mathbb{E}\bigl(\widehat
{f}%
(t_{0})-f(t_{0})\bigr)^{2}\leq C\cdot\biggl(\frac{\log n}{n}\biggr)^{({2\alpha})/({%
1+2\alpha})}.
\end{equation}
\end{theorem}

Theorem \ref{LocalAdapt.thm} shows that both estimators are spatially
adaptive, without prior knowledge of the smoothness of the underlying
functions.

\subsection{Regression in general natural exponential families}
\label{General-NEF.sec}

We have so far focused on the nonparametric regression in the NEF--QVF
families. Our method can be extended to the nonparametric regression in the
general one-parameter natural exponential families where the variance
is no
longer a quadratic function of the mean.

Suppose we observe
%
\begin{equation}\label{data}
Y_{i}\stackrel{\mathrm{ind}}{\sim}\operatorname{NEF}(f(t_{i})),\qquad
i=1,\ldots,n, t_{i}={\frac{i}{n}},
\end{equation}
and wish to estimate the mean function $f(t)$. When the variance is not a
quadratic function of the mean, the VST still exists, although the
mean-matching VST does not. In this case, we set $a=b=0$ in (\ref{Hm}) and
define $H_{m}$ as
%
\begin{equation}
H_{m}(X)=G\biggl({\frac{X}{m}}\biggr).
\end{equation}
We then apply the same four-step procedure, Binning--VST--Gaussian
Regres\-sion--Inverse VST, as outlined in Section \ref{procedure.sec} where
either BlockJS or NeighCoeff is used in the third step. Denote the resulting
estimator by $\widehat{f}_{\mathrm{BJS}}$ and $\widehat{f}_{\mathrm{NC}}$, respectively.

The following theorem is an extension of Theorem \ref{Yi.thm} to the general
one-parameter natural exponential families where the standard VST is used.
\begin{theorem}
\label{YiG.thm} Let $f\in F^{d}(M,\varepsilon,v)$. Then $Y_{j}^{\ast
}=G(%
\frac{Q_{j}}{m})$ can be written as
%
\begin{equation} \label{YiG.decomp}
Y_{j}^{\ast}=G\biggl(f\biggl({\frac{j}{T}}\biggr)\biggr)+\epsilon_{j}+m^{-{
{1/2}}}Z_{j}+\xi
_{j},\qquad j=1,2,\ldots,T,
\end{equation}
where $Z_{j}\stackrel{\mathrm{i.i.d.}}{\sim}N(0,1)$, $\epsilon_{j}$ are constants
satisfying $|\epsilon_{j}|\leq{c} ( m^{-1}+T^{-d} ) $ and
consequently for some constant $C>0$
%
\begin{equation} \label{epsilonG.bnd}
{\frac{1}{T}}\sum_{j=1}^{T}\epsilon_{j}^{2}\leq C (
m^{-2}+T^{-2d} )
\end{equation}
and $\xi_{j}$ are independent and ``stochastically
small'' random variables satisfying that for any integer $%
k>0$ and any constant $a>0$
%
\begin{eqnarray} \label{momentG.bd}
\mathbb{E}|\xi_{j}|^{k} &\leq& C_{k}\log^{2k}m\cdot(m^{-k}+T^{-dk}) %
\quad\mbox{and}\nonumber\\[-8pt]\\[-8pt]
\mathbb{P}(|\xi_{j}|>a) &\leq& C_{k}\log^{2k}m\cdot
(m^{-k}+T^{-dk})a^{-k},\nonumber
\end{eqnarray}
where $C_{k}>0$ is a constant depending only on $k,d$ and $M$.
\end{theorem}

The proof of Theorem \ref{YiG.thm} is similar to that of Theorem \ref
{Yi.thm}%
. Note that the bound for the deterministic error in (\ref
{epsilonG.bnd}) is
different from the one given in equation (\ref{epsilon.bnd}). This
difference affects the choice of the bin size.
\begin{theorem}
\label{NEF-global.thm} Suppose the wavelet $\psi$ is r-regular. Let $%
X_{i}\sim \operatorname{NEF}(f(t_{i})), i=1,\ldots,n, t_{i}={\frac{i}{n}}$. Let
$T=cn^{{1/2}}$. Then the estimator $\widehat{f}_{\mathrm{BJS}}$ satisfies
\[
\sup_{f\in F_{p,q}^{\alpha}(M,\varepsilon,v)}\mathbb{E}\Vert\widehat
{f}%
_{\mathrm{BJS}}-f\Vert_{2}^{2}\leq\cases{
Cn^{-{({2\alpha})/({1+2\alpha})}}, \cr
\qquad\mbox{$p\geq2, \alpha\leq r$
and $\biggl(\alpha-{\dfrac{1}{p}}\biggr)>{\dfrac{2\alpha}{1+2\alpha}}$},
\vspace*{2pt}\cr
Cn^{-{({2\alpha})/({1+2\alpha})}}(\log n)^{({2-p})/({p(1+2\alpha
)})}, \cr
\qquad\mbox{$1\leq p<2, \alpha\leq r$ and $\biggl(\alpha-{\dfrac{1}{p}}\biggr)>{\dfrac
{2\alpha}{1+2\alpha}}$},}
\]
and the estimator $\widehat{f}_{\mathrm{NC}}$ satisfies
\begin{eqnarray}
\sup_{f\in F_{p,q}^{\alpha}(M,\varepsilon,v)}\mathbb{E}\Vert\widehat
{f}%
_{\mathrm{NC}}-f\Vert_{2}^{2}\leq C \biggl( \frac{\log n}{n} \biggr) ^{{({2\alpha})/({%
1+2\alpha})}},\nonumber\\
\eqntext{p\geq1, \alpha\leq r \mbox{ and }{\biggl(\alpha-{\dfrac{1}{%
p}}\biggr)>{\dfrac{2\alpha}{1+2\alpha}}}.}
\end{eqnarray}
\end{theorem}
\begin{remark}
Note that the number of bins here is $T=O(n^{{{1/2}}})$.
This gives a larger bin size than that needed with NEF--QVF. Because
the VST yields higher bias than the
mean-matching VST in the case of NEF--QVF, it is necessary to use larger
bins. The condition $(\alpha-{\frac{1}{p}})>{\frac{2\alpha
}{1+2\alpha}} $ is
also stronger than the condition $\frac{3}{2}(\alpha-{\frac
{1}{p}})>{\frac{%
2\alpha}{1+2\alpha}}$ which is needed in the case of NEF--QVF. The
functions are required to be smoother than before. This is due to the fact
that both the approximation error and the discretization error are
larger in
this case. See Section \ref{discussion.sec} for more discussions.
\end{remark}

We have the following result on spatial adaptivity.
\begin{theorem}
\label{NEF-LocalAdapt.thm} Suppose the wavelet $\psi$ is r-regular
with ${%
\frac{1}{2}}<\alpha\leq r$. Let $t_{0}\in(0,1)$ be fixed. Let
$X_{i}\sim
\operatorname{NEF}(f(t_{i})), i=1,\ldots,n, t_{i}=\frac{i}{n}$. Let $T=cn^{{{1/2}}}$.
Then for $\widehat{f}=\widehat{f}_{\mathrm{BJS}}$ or $\widehat{f}_{\mathrm{NC}}$
%
\begin{equation}\label{NEF-point}
\sup_{f\in F^{\alpha}(M,t_{0},\delta,\varepsilon,v)}\mathbb
{E}\bigl(\widehat{f}%
(t_{0})-f(t_{0})\bigr)^{2}\leq C\biggl(\frac{\log n}{n}\biggr)^{({2\alpha
})/({1+2\alpha})}.
\end{equation}
\end{theorem}
\begin{remark}
In Remark \ref{remark1} we noted that some nonexponential families
admit mean-matching
variance stabilizing transformations. Although we do not pursue the
issue in the
current paper, we believe that analogs of our procedure can be developed
for these families and the basic results in Theorems \ref{global.thm}
and \ref{LocalAdapt.thm} can be extended to such
situations. A different possibility is that the error distributions lie
in a one parameter
family that admits a VST that is not mean matching. In that case one
could expect
analogs of Theorems \ref{NEF-global.thm} and \ref{NEF-LocalAdapt.thm}
to be valid.
\end{remark}

\subsection{Discussion}
\label{discussion.sec}

Our procedure begins with binning. This step makes the data more ``normal''
and at the same time reduces the number of observations from $n$ to $T$.
This step in general does not affect the rate of convergence as long as the
underlying function has certain minimum smoothness so that the bias induced
by local averaging is negligible relative to the minimax estimation risk.
While the number of observations is reduced by binning, the noise level is
also reduced accordingly.

An important quantity in our method is the value of $T$, the number of bins,
or equivalently the value of the bin size $m$. The choice of $T= cn^{3/4}$
for the NEF--QVF and $T=cn^{1/2}$ for the general NEF are determined by the
bounds for the approximation error, the discretization\vspace*{1pt} error, and the
stochastic error. For functions in the Besov ball $B^\alpha_{p,q}(M)$, the
discretization error between the sampled function $\{G(f({j/T}))\dvtx j=1,\ldots,
T\}$ and the whole\vspace*{1pt} function $G(f(t))$ can be bounded by $C T^{-2d}$
where $%
d=(\alpha-{\frac{1}{p}})\wedge1$ (see\vspace*{-2pt} Lemma \ref{besov.approx.lem} in
Section \ref{SingleBlockRisk.sec}). The approximation error ${\frac
{1}{T}}%
\sum_{i=1}^T\epsilon_i^2$ can be bounded by $C(m^{-4} + T^{-2d})$ as
in (\ref%
{epsilon.bnd}). In order to adaptively achieve the optimal rate of
convergence, these deterministic errors need to be negligible relative to
the minimax rate of convergence $n^{-{({2\alpha})/({1+2\alpha})}}$
for all $%
\alpha$ under consideration. That is, we need to have
$m^{-4}=o(n^{-{({%
2\alpha})/({1+2\alpha})}})$ and $T^{-2d}=o(n^{-{({2\alpha
})/({1+2\alpha})}})$. These conditions put constraints on both $m $ and $\alpha$ (and
$p$). We
choose $m=cn^{{1/4}}$ (or equivalently $T=cn^{{3/4}}$) to ensure
that the approximation error is always negligible for all $\alpha$. This
choice also guarantees that the stochastic error is under control. With this
choice of $m$, we then need $\frac{3}{2}(\alpha-{\frac
{1}{p}})>{\frac{%
2\alpha}{1+2\alpha}}$ or equivalently $\frac{2\alpha^{2}-\alpha/3}{
1+2\alpha}>{\frac{1}{p}}$.

In the natural exponential family with a quadratic variance function, the
existence of a mean-matching VST makes the approximation error small and
this provides advantage over more general natural exponential\vspace*{1pt} families. For
general NEF without a quadratic variance function, the approximation
error ${%
\frac{1}{T}}\sum_{i=1}^T\epsilon_i^2$ is of order $m^{-2}+ T^{-2d}$ instead
of $m^{-4} + T^{-2d}$. Making it negligible for all $\alpha$ under
consideration requires $m=cn^{{{1/2}}}$. With this\vspace*{-2pt} choice of
$m$, we
require $\alpha-{\frac{1}{p}}>{\frac{2\alpha}{1+2\alpha}}$ or equivalently
$\frac{2\alpha^{2}-\alpha}{1+2\alpha}>{\frac{1}{p}}$ in order to control
the discretization error. In particular, this condition is satisfied if
$%
\alpha\ge1 + {\frac{1}{p}}$.

In this paper we present a unified approach to nonparametric regression in
the natural exponential families and the optimality results are given for
Besov spaces. As mentioned in the \hyperref[intro]{Introduction}, a
wavelet shrinkage and
modulation method was introduced in Antoniadis and Sapatinas (\citeyear{AntSap01}) for
regression in the NEF--QVF and it was shown that the estimator attains the
optimal rate over the classical Sobolev spaces with the smoothness
index $%
\alpha> 1/2$. In comparison to the results given in Antoniadis and
Sapatinas (\citeyear{AntSap01}), our results are more general in terms of the function
spaces as well as the natural exponential families. On the other hand, we
require slightly stronger conditions on the smoothness of the underlying
functions. It is\vspace*{1pt} intuitively clear that through binning and VST a certain
amount of bias is introduced. The conditions\vspace*{-1pt} $\frac{3}{2}(\alpha
-{\frac{1}{p%
}})>{\frac{2\alpha}{1+2\alpha}}$ in the case of NEF--QVF and $\alpha
-{%
\frac{1}{p}}>{\frac{2\alpha}{1+2\alpha}}$ in the general case are the
minimum smoothness condition needed to ensure that the bias is under
control. The bias in the general NEF case is larger and therefore the
required smoothness condition is stronger.

\section{Numerical study}
\label{numerical.sec}

In this section we study the numerical performance of our estimators. The
procedures introduced in Section \ref{procedure.sec} are easily
implementable. We shall first consider simulation results and then
apply one
of our procedures in the analysis of two real data sets.

\subsection{Simulation results}

As discussed the Section \ref{vst.sec}, there are several different versions
of the VST in the literature and we have emphasized the importance of using
the mean-matching VST for theoretical reasons. We shall now consider the
effect of the choice of the VST on the numerical performance of the
resulting estimator. To save space we only consider the Poisson and
Bernoulli cases. We shall compare the numerical performance of the
mean-matching VST with those of classical transformations by Bartlett (\citeyear{Bartlett36})
and Anscombe (\citeyear{Anscombe48}) using simulations. The transformation formulae are given
as follows. (In the following tables and figures, we shall use MM for
mean-matching.)\vspace*{-10pt}

\begin{table}[h]
\begin{tabular*}{\tablewidth}{@{\extracolsep{\fill}}lccc@{}}
\hline
& \textbf{MM} & \textbf{Bartlett} & \textbf{Anscombe} \\
\hline
$\operatorname{Poi}(\lambda) $ & $\sqrt{X+1/4}$ & $\sqrt{X}$ & $\sqrt
{X+3/8}$ \\[2pt]
$\operatorname{Bin}(m,p) $ & $\sin^{-1}\sqrt{\frac{X+1/4}{m+1/2}}$ &
$\sin
^{-1}\sqrt{\frac{X}{m}}$ & $\sin^{-1}\sqrt{\frac{X+3/8}{m+3/4}}$ \\
\hline
\end{tabular*}\vspace*{-10pt}
\end{table}

\begin{figure}

\includegraphics{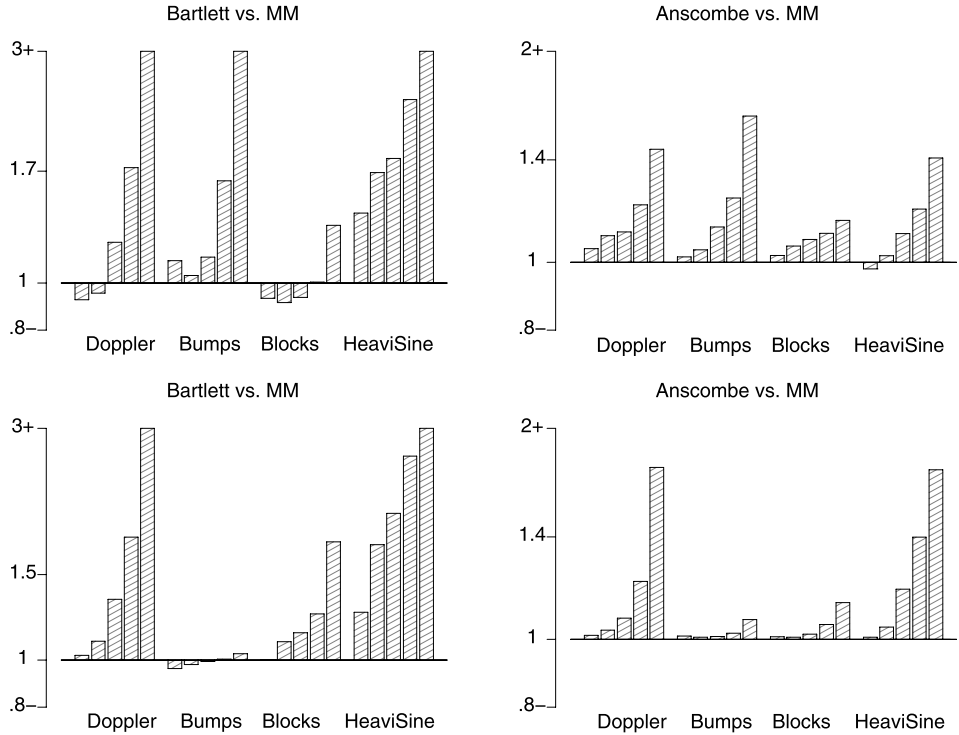}

\caption{Left panels: the vertical bars represent the ratios
of the MSE of the estimator using the Bartlett VST to the corresponding MSE
of our estimator using the mean-matching VST. Right Panels: the bars
represent the ratios of the MSE of the estimator using the Anscombe VST to
the corresponding MSE of the estimator using the mean-matching VST. The
higher the bar the better the relative performance of our estimator. The
bars are plotted on a log scale and the original ratios are truncated
at the
value $3$ for the Bartlett VST and at $2$ for the Anscombe VST. For
each signal
the bars are ordered from left to right in the order of increasing sample
size. The top row is for the Bernoulli case and the bottom row for the
Poisson case.}
\label{VST-risk-comp}
\end{figure}

Four standard test functions, Doppler, Bumps, Blocks and HeaviSine,
representing different level of spatial variability are used for the
comparison of the three VSTs. See Donoho and Johnstone (\citeyear{DJ94}) for the
formulae of the four test functions. These test functions are suitably
normalized so that they are positive and taking values between 0 and 1 (in
the binomial case). Sample sizes vary from a few hundred to a few hundred
thousand. We use Daubechies' compactly supported wavelet \textit{Symmlet} 8 for
wavelet transformation. As is the case in general, it is possible to obtain
better estimates with different wavelets for different signals. But for
uniformity, we use the same wavelet for all cases. Although our asymptotic
theory only gives a justification for the choice of the bin size of
order $%
n^{1/4}$ due to technical reasons, our extensive numerical studies have
shown that the procedure works well when the number of counts in each
bin is
between 5 and 10 for the Poisson case, and similarly for the Bernoulli case
the average number of successes and failures in each bin is between 5 and
10. We follow this guideline in our simulation study. Table \ref{simu.table}
reports the average squared errors over 100 replications for the BlockJS
thresholding. The sample sizes are $1280, 5120, \ldots, 327\mbox{,}680$ for the
Bernoulli
case and $640, 2560, \ldots, 163\mbox{,}840$ for the Poisson case. A graphical
presentation is given in Figure \ref{VST-risk-comp}.

\begin{table}
\caption{Mean squared error (MSE) from $100$ replications. The MSE is
in units of $10^{-3}$ for~Bernoulli~case~and~$10^{-2}$ for Poisson case}
\label{simu.table}
\begin{tabular*}{\tablewidth}{@{\extracolsep{\fill}}ld{2.3}d{2.3}d{2.3}cd{3.3}d{3.3}d{3.3}@{}}
\hline
& \multicolumn{1}{c}{\textbf{MM}} & \multicolumn{1}{c}{\textbf{Bartlett}} &
\multicolumn{1}{c}{\textbf{Anscombe}}
& & \multicolumn{1}{c}{\textbf{MM}} & \multicolumn{1}{c}{\textbf{Bartlett}} & \multicolumn{1}{c@{}}{\textbf{Anscombe}}\\
\hline
\multicolumn{8}{c}{Bernoulli}\\[4pt]
Doppler & & & & Bumps \\[2pt]
\phantom{00,}1280 & 12.117 & 11.197 & 12.673 & \phantom{00,}1280 & 7.756 & 8.631 & 7.896 \\
\phantom{00,}5120 & 3.767 & 3.593 & 4.110 & \phantom{00,}5120 & 7.455 & 7.733 & 7.768 \\
\phantom{0}20,480 & 1.282 & 1.556 & 1.417 & \phantom{0}20,480 & 3.073 & 3.476 & 3.450 \\
\phantom{0}81,920 & 0.447 & 0.772 & 0.540 & \phantom{0}81,920 & 1.203 & 1.953 & 1.485 \\
327,680 & 0.116 & 0.528 & 0.169 & 327,680 & 0.331 & 1.312 & 0.535\\[2pt]
Blocks & & & & HeaviSine\\[2pt]
\phantom{00,}1280 & 18.451 & 17.171 & 18.875 & \phantom{00,}1280 & 2.129 & 2.966 & 2.083 \\
\phantom{00,}5120 & 7.582 & 6.911 & 7.996 & \phantom{00,}5120 & 0.842 & 1.422 & 0.860 \\
\phantom{0}20,480 & 3.288 & 3.072 & 3.545 & \phantom{0}20,480 & 0.549 & 0.992 & 0.603 \\
\phantom{0}81,920 & 1.580 & 1.587 & 1.737 & \phantom{0}81,920 & 0.285 & 0.681 & 0.339 \\
327,680 & 0.594 & 0.781 & 0.681 & 327,680 & 0.138 & 0.532 & 0.195\\[4pt]
\multicolumn{8}{c}{Poisson}\\[4pt]
Doppler & & & & Bumps\\[2pt]
\phantom{000,}640 & 8.101 & 8.282 & 8.205 & \phantom{000,}640 & 107.860 & 103.696 & 109.023 \\
\phantom{00,}2560 & 3.066 & 3.352 & 3.160 & \phantom{00,}2560 & 70.034 & 68.616 & 70.495 \\
\phantom{0}10,240 & 1.069 & 1.426 & 1.146 & \phantom{0}10,240 & 24.427 & 24.268 & 24.653 \\
\phantom{0}40,960 & 0.415 & 0.743 & 0.502 & \phantom{0}40,960 & 9.427 & 9.469 & 9.620 \\
163,840 & 0.108 & 0.461 & 0.190 & 163840 & 3.004 & 3.098 & 3.204\\[2pt]
Blocks & & & & HeaviSine\\[2pt]
\phantom{000,}640 & 12.219 & 12.250 & 12.320 & \phantom{000,}640 & 2.831 & 3.552 & 2.851 \\
\phantom{00,}2560 & 5.687 & 6.209 & 5.724 & \phantom{00,}2560 & 0.849 & 1.468 & 0.884 \\
\phantom{0}10,240 & 2.955 & 3.363 & 3.005 & \phantom{0}10,240 & 0.425 & 0.852 & 0.501 \\
\phantom{0}40,960 & 1.424 & 1.773 & 1.495 & \phantom{0}40,960 & 0.213 & 0.560 & 0.298 \\
163,840 & 0.508 & 0.890 & 0.573 & 163,840 & 0.118 & 0.455 & 0.206\\
\hline
\end{tabular*}
\end{table}

Table \ref{simu.table} compares the performance of three nonparametric
function estimators constructed from three VSTs and wavelet BlockJS
thresholding for Bernoulli and Poisson regressions. The three VSTs are the
mean-matching, Bartlett and Anscombe transformations given above. The
results show the mean-matching VST outperforms the classical transformations
for nonparametric estimation in most cases. The improvement becomes more
significant as the sample size increases.

In the Poisson regression, the mean-matching VST outperforms the Bartlett
VST in 17 out of 20 cases and uniformly outperforms the Anscombe VST in all
20 cases. The case of Bernoulli regression is similar: the
mean-matching VST
is better than the Bartlett VST in 15 out of 20 cases and better than the
Anscombe VST in 19 out of 20 cases. Although the mean-matching VST does not
uniformly dominate either the Bartlett VST or the Anscombe VST, the
improvement of the mean-matching VST over the other two VSTs is significant
as the sample size increases for all four test functions. The simulation
results show that mean-matching VST yields good numerical results in
comparison to other VSTs. These numerical findings is consistent with the
theoretical results given in Section \ref{theory.sec} which show that the
estimator constructed from the mean-matching VST enjoys desirable adaptivity
properties.

Table \ref{simu2.table} reports the average squared errors over 100
replications for the NeighCoeff procedure in the same setting as those in
Table \ref{simu.table}. In comparison to BlockJS, the numerical performance
of NeighCoeff is overall slightly better. Among the three VSTs, the
mean-matching VST again outperforms both the Anscombe VST and Bartlett VST.

\begin{figure}

\includegraphics{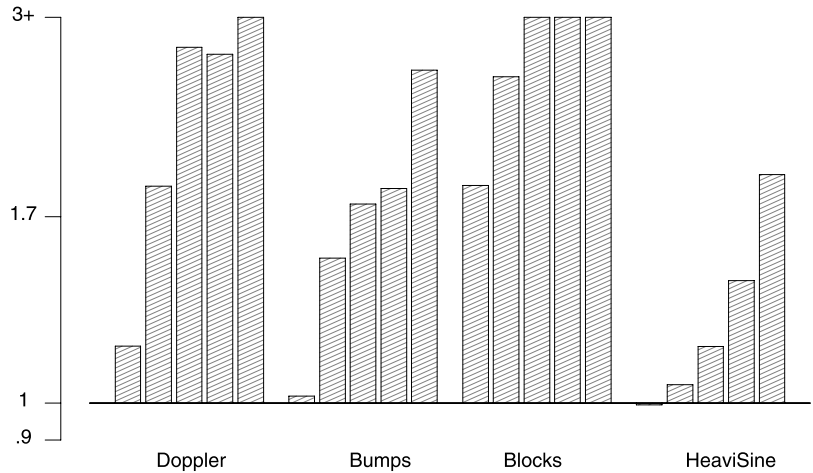}

\caption{The vertical bars represent the ratios of the MSE
of Donoho's estimator to the corresponding MSE of our estimator. The higher
the bar the better the relative performance of our estimator. The bars are
plotted on a log scale and the original ratios are truncated at the
value $3$.
For each signal the bars are ordered from left to right in the order of
increasing sample size.}
\label{VST-risk-comp2}
\end{figure}

We have so far considered the effect of the choice of VST on the performance
of the estimator. We now discuss the Poisson case in more detail and compare
the numerical performance of our procedure with other estimators
proposed in
the literature. As mentioned in the \hyperref[intro]{Introduction},
Besbeas, De Feis and
Sapatinas (\citeyear{BFS04}) carried out an extensive simulation studies comparing
several nonparametric Poisson regression estimators including the estimator
given in Donoho (\citeyear{Donoho93}). The estimator in Donoho (\citeyear{Donoho93}) was constructed by
first applying the Anscombe (\citeyear{Anscombe48}) VST to the binned data and by then using
a wavelet procedure with a global threshold such as VisuShrink [Donoho and
Johnstone (\citeyear{DJ94})] to the transformed data as if the data were actually
Gaussian. Figure \ref{VST-risk-comp2} plots the ratios of the MSE of
Donoho's estimator to the corresponding MSE of our estimator. The results
show that our estimator outperforms Donoho's estimator in all but one case
and in many cases our estimator has the MSE less than one half and sometimes
even one third of that of Donoho's estimator.

Besbeas, De Feis and Sapatinas (\citeyear{BFS04}) plotted simulation results of 27
procedures for six intensity functions (Smooth, Angles, Clipped Blocks,
Bumps, Spikes and Bursts) with sample size 512 under the squared root of
mean squared error (RMSE). We apply NeighCoeff and BlockJS procedures to
data with exactly the same intensity functions. The following table reports
the RMSE of NeighCoeff and BlockJS procedures based on 100 replications:  

We compare our results with the plots of RMSE for 27 methods in
Besbeas, De
Feis and Sapatinas (\citeyear{BFS04}). The NeighCoeff procedure dominates all 27 methods
for signals Smooth and Spikes, outperforms most of procedures for signals
Angles and Bursts, and performs slightly worse than average for signals
Clipped Blocks and Bumps. The BlockJS procedure is comparable with the
NeighCoeff procedure except for two signals Clipped Blocks and Bumps. We
should note that an exact numerical comparison here is difficult as the
results in Besbeas, de Feis and Sapatinas (\citeyear{BFS04}) were given in plots, not
numerical values.

\begin{table}
\caption{Mean squared error (MSE) from $100$ replications for the NeighCoeff
thresholding. The MSE is in units of $10^{-3}$ for Bernoulli case and $%
10^{-2}$ for Poisson case}
\label{simu2.table}
\begin{tabular*}{\tablewidth}{@{\extracolsep{\fill}}ld{2.3}d{2.3}d{2.3}cd{3.3}d{3.3}d{3.3}@{}}
\hline
& \multicolumn{1}{c}{\textbf{MM}} & \multicolumn{1}{c}{\textbf{Bartlett}} &
\multicolumn{1}{c}{\textbf{Anscombe}}
& & \multicolumn{1}{c}{\textbf{MM}} & \multicolumn{1}{c}{\textbf{Bartlett}} & \multicolumn{1}{c@{}}{\textbf{Anscombe}}\\
\hline
\multicolumn{8}{c}{Bernoulli}\\[4pt]
Doppler & & & & Bumps \\[2pt]
\phantom{00,}1280 & 8.574 & 8.569 & 8.959 & \phantom{00,}1280 & 7.085 & 7.741 & 7.361 \\
\phantom{00,}5120 & 2.935 & 3.211 & 3.129 & \phantom{00,}5120 & 6.810 & 7.052 & 7.180 \\
\phantom{0}20,480 & 1.029 & 1.380 & 1.143 & \phantom{0}20,480 & 2.846 & 3.364 & 3.204 \\
\phantom{0}81,920 & 0.377 & 0.800 & 0.438 & \phantom{0}81,920 & 0.958 & 1.789 & 1.220 \\
327,680 & 0.138 & 0.556 & 0.186 & 327,680 & 0.264 & 1.274 & 0.458\\[2pt]
Blocks & & & & HeaviSine\\[2pt]
\phantom{00,}1280 & 14.838 & 13.964 & 15.336 & \phantom{00,}1280 & 2.072 & 3.092 & 2.010 \\
\phantom{00,}5120 & 7.129 & 6.615 & 7.511 & \phantom{00,}5120 & 0.822 & 1.479 & 0.841 \\
\phantom{0}20,480 & 3.131 & 2.904 & 3.388 & \phantom{0}20,480 & 0.529 & 1.007 & 0.580 \\
\phantom{0}81,920 & 1.266 & 1.350 & 1.400 & \phantom{0}81,920 & 0.235 & 0.660 & 0.286 \\
327,680 & 0.469 & 0.680 & 0.553 & 327,680 & 0.102 & 0.512 & 0.156
\\[4pt]
\multicolumn{8}{c}{Poisson}\\[4pt]
Doppler & & & & Bumps\\[2pt]
\phantom{000,}640 & 7.789 & 8.030 & 7.888 & \phantom{000,}640 & 105.624 & 101.486 & 106.76 \\
\phantom{00,}2560 & 3.112 & 3.398 & 3.200 & \phantom{00,}2560 & 69.627 & 68.175 & 70.105 \\
\phantom{0}10,240 & 1.006 & 1.362 & 1.081 & \phantom{0}10,240 & 24.448 & 24.304 & 24.672 \\
\phantom{0}40,960 & 0.402 & 0.731 & 0.488 & \phantom{0}40,960 & 9.312 & 9.341 & 9.507 \\
163,840 & 0.106 & 0.460 & 0.187 & 163,840 & 3.005 & 3.102 & 3.203\\
\\[2pt]
Blocks & & & & HeaviSine\\[2pt]
\phantom{000,}640 & 12.301 & 12.141 & 12.412 & \phantom{000,}640 & 2.679 & 3.465 & 2.672 \\
\phantom{00,}2560 & 5.719 & 6.229 & 5.758 & \phantom{00,}2560 & 0.903 & 1.427 & 0.977 \\
\phantom{0}10,240 & 2.985 & 3.363 & 3.046 & \phantom{0}10,240 & 0.429 & 0.852 & 0.505 \\
\phantom{0}40,960 & 1.399 & 1.755 & 1.469 & \phantom{0}40,960 & 0.215 & 0.562 & 0.300 \\
163,840 & 0.504 & 0.877 & 0.572 & 163,840 & 0.120 & 0.453 & 0.209\\
\hline
\end{tabular*}
\end{table}
\begin{table} 
\begin{tabular*}{\tablewidth}{@{\extracolsep{\fill}}lcccccc@{}}
\hline
& \textbf{Smooth} & \textbf{Angles} & \textbf{Clipped blocks} & \textbf{Bumps}
& \textbf{Spikes} & \textbf{Bursts} \\
\hline
NeighCoeff & 1.773 & 2.249 & 5.651 & 4.653 & 2.096 & 2.591 \\
BlockJS & 1.760 & 2.240 & 6.492 & 5.454 & 2.315 & 2.853\\
\hline
\end{tabular*}\vspace*{-10pt}
\end{table}

\subsection{Real data applications}

We now demonstrate our estimation method in the analysis of two real data
sets, a gamma-ray burst data set (GRBs) and a packet loss data set. These
two data sets have been discussed in Kolaczyk and Nowak (\citeyear{KN05}).

%
\begin{figure}[b]

\includegraphics{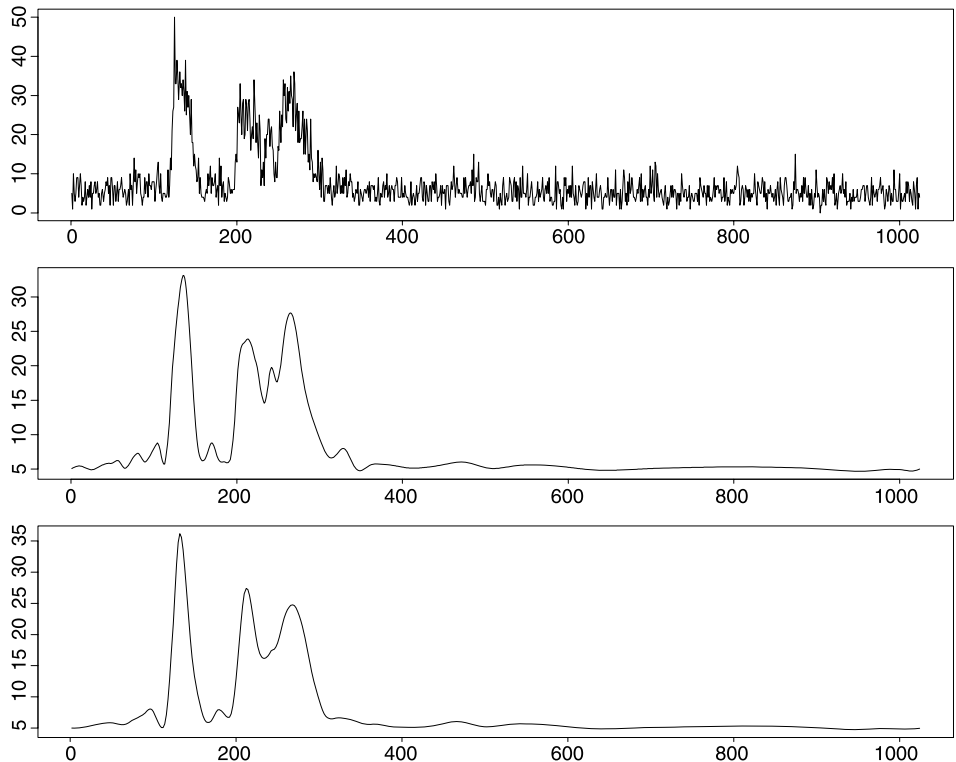}

\caption{Gamma-ray burst. Histogram of BATSE $551$ with $1024$
bins (top panel). Estimator based on $1024$ bin (middle panel).
Estimator with
$512$ bins (bottom panel).}
\label{gamma-ray-fig}
\end{figure}

Cosmic gamma-ray bursts were first discovered in the late 1960s. In 1991,
NASA launched the Compton Gamma Ray Observatory and its Burst and Transient
Source Explorer (BATSE) instrument, a sensitive gamma-ray detector. Much
burst data has been collected since then, followed by extensive studies and
many important scientific discoveries during the past few decades; however,
the source of GRBs remains unknown [Kaneko (\citeyear{Kaneko05})]. For more details
see the
NASA website \href{http://www.batse.msfc.nasa.gov/batse/}{http://www.batse.msfc.nasa.gov/batse/}. GRBs seem to be
connected to massive stars and become powerful probes of the star formation
history of the universe. However not many redshifts are known and there is
still much work to be done to determine the mechanisms that produce these
enigmatic events. Statistical methods for temporal studies are
necessary to
characterize their properties and hence to identify the physical properties
of the emission mechanism. One of the difficulties in analyzing the time
profiles of GRBs is the transient nature of GRBs which means that the usual
assumptions for Fourier transform techniques do not hold [Quilligan et al.
(\citeyear{Quilliganetal02})]. We may model the time series data by an inhomogeneous Poisson
process, and apply our wavelet procedure. The data set we use is called
BATSE 551 with the sample size 7808. In Figure \ref{gamma-ray-fig},
the top
panel is the histogram of the data with 1024 bins such that the number of
observations in each bin would be between 5 and 10. In fact we have on
average 7.6 observations. The middle panel is the estimate of the intensity
function using our procedure. If we double the width of each bin, that
is, the
total number of bins is now 512, the new estimator in the bottom panel is
noticeably different from previous one since it does not capture the fine
structure from time 200 to 300. The study of the number of pulses in GRBs
and their time structure is important to provide evidence for rotation
powered systems with intense magnetic fields and the added complexity
of a
jet.

\begin{figure}[b]

\includegraphics[scale=0.99]{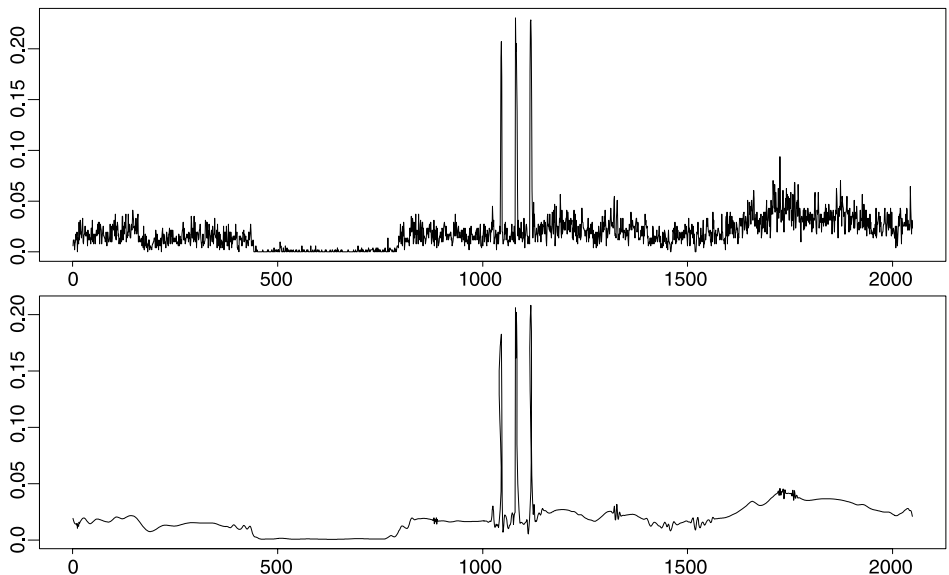}

\caption{Packet loss data. Histogram with $2048$ bins (top
panel). Estimator based on the binned data (bottom panel).}
\label{packet-loss}
\end{figure}

Packet loss describes an error condition in internet traffic in which data
packets appear to be transmitted correctly at one end of a connection, but
never arrive at the other. So, if 10 packets were sent out, but only 8 made
it through, then there would be 20\% overall packet loss. The following data
were originally collected and analyzed by Yajnik et al. (\citeyear{Yajniketal99}). The
objective is to understand packet loss by modeling. It measures the
reliability of a connection and is of fundamental importance in network
applications such as audio/video conferencing and Internet telephony.
Understanding the loss seen by such applications is important in their
design and performance analysis. The measurements are of loss as seen by
packet probes sent at regular time intervals. The packets were transmitted
from the University of Massachusetts at Amherst to the Swedish
Institute of
Computer Science. The records note whether each packet arrived or was lost.
It is a Bernoulli time series, and can be naturally modeled as Binomial
after binning the data. Figure \ref{packet-loss} gives the histogram and our
corresponding estimator. The average sum of failures in each bin is about
10. The estimator in Kolaczyk and Nowak (\citeyear{KN05}) is comparable to ours. But
our procedure is more easily implemented.

\section{Proofs}
\label{proof.sec}

In this section we give proofs for Theorems \ref{Yi.thm}, \ref{global.thm}
and \ref{NEF-global.thm}. Theorems~\ref{LocalAdapt.thm} and \ref%
{NEF-LocalAdapt.thm} can be proved in a similar way as Theorem 4 in Brown,
Cai and Zhou (\citeyear{BCZ08}) by applying Proposition \ref{single.block.prop} in
Section \ref{SingleBlockRisk.sec}. We begin by proving Lemmas \ref%
{approx.lem} and \ref{xi.lem} as well as an additional technical result,
Lemma \ref{iidapprox}. These results are needed to establish Theorem
\ref%
{Yi.thm} in which an approximation bound between our model and a Gaussian
regression model is given explicitly. Finally we apply Theorem \ref{Yi.thm}
and risk bounds for block thresholding estimators in Proposition \ref%
{single.block.prop} to prove Theorems~\ref{global.thm} and \ref%
{NEF-global.thm}.

\subsection{Proof of preparatory technical results}
\label{coupling-proof.sec}

\mbox{}

\begin{pf*}{Proof of Lemma \protect\ref{approx.lem}}
We only prove (\ref{bias}), the first part of the lemma. The
proof for
equation (\ref{var}), the second part, is similar and simpler. By Taylor's
expansion we write
\[
G \biggl( \frac{X+a}{m+b} \biggr) -G ( \mu(\eta) )
=T_{1}+T_{2}+T_{3}+T_{4},
\]
where
\begin{eqnarray*}
T_{1}&=&G^{\prime} ( \mu(\eta) ) \biggl( \frac{X+a}{m+b}-\mu(\eta
) \biggr) ,\qquad T_{2}=\frac{1}{2}G^{\prime\prime} ( \mu(\eta
) ) \biggl( \frac{X+a}{m+b}-\mu(\eta) \biggr) ^{2}, \\
T_{3} &=& \frac{1}{6}G^{\prime\prime\prime} ( \mu(\eta) ) \biggl(
\frac{X+a}{m+b}-\mu(\eta) \biggr) ^{3}, \qquad T_{4}=\frac{1}{24}%
G^{(4)} ( \mu^{\ast} ) \biggl( \frac{X+a}{m+b}-\mu(\eta) \biggr)
^{4}
\end{eqnarray*}
and $\mu^{\ast}$ is in between $\frac{X+a}{m+b}$ and $\mu(\eta)$. By
definition, $G^{\prime} ( \mu(\eta) ) =I ( \eta)
^{-1/2}$ with $I ( \eta) =\mu^{\prime} ( \eta) $
which is also $V(\mu( \eta) )$ in (\ref{vst}), then
\[
G^{\prime\prime} ( \mu(\eta) ) \mu^{\prime} ( \eta
) =-\tfrac{1}{2}I ( \eta) ^{-3/2}I^{\prime} ( \eta
),
\]
that is,
\[
G^{\prime\prime} ( \mu(\eta) ) =-\tfrac{1}{2}I ( \eta
) ^{-5/2}I^{\prime} ( \eta),
\]
then%
\begin{eqnarray*}
\mathbb{E}T_{1} &=&I ( \eta) ^{-1/2}\frac{a-\mu(\eta)b}{m+b}, \\
\mathbb{E}T_{2} &=&-\frac{1}{4}I ( \eta) ^{-5/2}I'
( \eta) \biggl[ \biggl( \frac{a-\mu(\eta)b}{m+b} \biggr) ^{2}+%
\frac{mI ( \eta) }{ ( m+b ) ^{2}} \biggr].
\end{eqnarray*}
Note that $G''( \mu(\eta) ) $ is uniformly bounded on $\Theta$ by the
assumption in the lemma, then we have
%
\begin{eqnarray} \label{diffeq}
\mathbb{E} ( T_{1}+T_{2} ) &=&\frac{m}{ ( m+b )
^{2}I ( \eta) ^{1/2}} \biggl( a-\mu(\eta)b-\frac{\mu^{\prime
\prime} ( \eta) }{4\mu^{\prime} ( \eta) } \biggr)
+O \biggl( \frac{1}{m^{2}} \biggr) \nonumber\\[-8pt]\\[-8pt]
&=&\frac{1}{mI ( \eta) ^{1/2}} \biggl( a-\mu(\eta)b-\frac{\mu
^{\prime\prime} ( \eta) }{4\mu^{\prime} ( \eta) }%
\biggr) +O \biggl( \frac{1}{m^{2}} \biggr).\nonumber
\end{eqnarray}
It is easy to show that
%
\begin{equation}\label{diffeq1}
\vert\mathbb{E}T_{3} \vert= \biggl\vert\frac{1}{6}G^{\prime
\prime\prime} ( \mu(\eta) ) \mathbb{E} \biggl( \frac{X+a}{m+b}%
-\mu(\eta) \biggr) ^{3} \biggr\vert=O \biggl( \frac{1}{m^{2}} \biggr) ,
\end{equation}
since $ \vert\mathbb{E} ( X/m-\mu(\eta) ) ^{3} \vert
=O ( \frac{1}{m^{2}} ) $. For any $\epsilon>0$ it is known that
\[
\mathbb{P} \biggl\{ \biggl\vert\frac{X+a}{m+b}-\mu( \eta)
\biggr\vert>\epsilon\biggr\} \leq\mathbb{P} \{ \vert X/m-\mu
( \eta) \vert>\epsilon/2 \},
\]
which decays
exponentially fast as $m\rightarrow\infty$ [see, e.g., Petrov (\citeyear{Petrov75})]. This
implies $\mu^{\ast}$ is in the interior of the natural parameter
space and
then $G^{(4)} ( \mu^{\ast} ) $ is bounded with probability
approaching to $1$ exponentially fast. Thus we have
%
\begin{equation}\label{diffeq2}
\vert\mathbb{E}T_{4} \vert\leq C\mathbb{E} \biggl( \frac{X+a}{m+b}%
-\mu(\eta) \biggr) ^{4}=O \biggl( \frac{1}{m^{2}} \biggr).
\end{equation}
Equation (\ref{bias}) then follows immediately by combining equations
(\ref%
{diffeq})--(\ref{diffeq2}).
\end{pf*}
\begin{pf*}{Proof of Lemma \protect\ref{KMT.lem}}
The proof is similar to Corollary 1 of Zhou (\citeyear{zhou06}). Let $\widetilde
{X}=\frac{X-m\mu
}{\sqrt{mV}}$. It is shown in Koml\'{o}s, Major and Tusn\'{a}dy (\citeyear{KMT75}) that
there exists a standard normal random variable $Z\sim N(0,1)$ and
constants $%
\varepsilon,c_{4}>0$ not depending on $m$ such that whenever the event
$%
A=\{|\widetilde{X}|\leq\varepsilon\sqrt{m}\}$ occurs,
%
\begin{equation} \label{KMT1}
|\widetilde{X}-Z|<\frac{c_{4}}{\sqrt{m}}+\frac{c_{4}}{\sqrt{m}}\widetilde
{X}^{2}.
\end{equation}
Obviously inequality (\ref{KMT1}) still holds when $ \vert\widetilde{X}%
\vert\leq\varepsilon_{1}\sqrt{m}$ for $0<\varepsilon_{1}\leq
\varepsilon$. Let's choose $\varepsilon_{1}$ small enough such that $%
c_{4}\varepsilon_{1}^{2}<1/2$. When $ \vert\widetilde{X} \vert\leq
\varepsilon_{1}\sqrt{m}$, we have $ \vert\widetilde{X}-Z \vert\leq
\frac{c_{4}}{\sqrt{m}}+\frac{1}{2} \vert\widetilde{X} \vert$ from
(\ref%
{KMT1}), which\vspace*{-2pt} implies $ \vert\widetilde{X} \vert- \vert
Z \vert\leq\frac{c_{4}}{\sqrt{m}}+\frac{1}{2} \vert\widetilde{X}%
\vert$ by the triangle inequality, that is, $ \vert\widetilde{X}%
\vert\leq\frac{2c_{4}}{\sqrt{m}}+2 \vert Z \vert$, so we
have
\[
\vert\widetilde{X}-Z \vert\leq\frac{c_{4}}{\sqrt{m}}+\frac{c_{4}}{%
\sqrt{m}} \biggl( \frac{2c_{4}}{\sqrt{m}}+2 \vert Z \vert\biggr)
^{2}\leq c_{2}Z^{2}+c_{3}
\]
for some constants $c_{1},c_{2}>0$. 
\end{pf*}
\begin{pf*}{Proof of Lemma \protect\ref{xi.lem}}
By Taylor's expansion we write
\[
G \biggl( \frac{X+a}{m+b} \biggr) -G ( \mu) =G^{\prime} ( \mu
) \biggl( \frac{X+a}{m+b}-\mu\biggr) +\frac{1}{2}G^{\prime\prime
} ( \mu^{\ast} ) \biggl( \frac{X+a}{m+b}-\mu\biggr) ^{2}.
\]
Recall that $ \vert\epsilon\vert= \vert\mathbb{E}G (
\frac{X+a}{m+b} ) -G({\mu}) \vert=O(m^{-2})$ from Lemma \ref%
{approx.lem}, and $Z$ is a standard normal variable satisfying (\ref{KMT}),
and
%
\begin{equation}
\xi=G\biggl({\frac{X+a}{m+b}}\biggr)-G({\mu})-\epsilon-m^{-{{1/2}}}Z.
\end{equation}
We write $\xi=\xi_{1}+\xi_{2}+\xi_{3}$, where
\begin{eqnarray*}
\xi_{1} &=&G^{\prime} ( \mu) \biggl( \frac{X+a}{m+b}-\frac{X}{m}%
\biggr) -\epsilon=G^{\prime} ( \mu) \frac{am-bX}{m(m+b)}%
-\epsilon, \\
\xi_{2} &=&G^{\prime} ( \mu) \Biggl( \frac{X}{m}-\mu-\sqrt{%
\frac{V}{m}}Z \Biggr) =\frac{G^{\prime} ( \mu) }{m} \bigl( X-m\mu-%
\sqrt{mV}Z \bigr), \\
\xi_{3} &=&\frac{1}{2}G^{\prime\prime} ( \mu^{\ast} ) \biggl(
\frac{X+a}{m+b}-\mu\biggr) ^{2}=\frac{1}{2}G^{\prime\prime} ( \mu
^{\ast} ) \biggl( \frac{X-m\mu}{m+b}+\frac{a-b\mu}{m+b} \biggr) ^{2}.
\end{eqnarray*}
It is easy to see that $\mathbb{E}|\xi_{1}|^{k}\leq C_{k}m^{-k}$.
Since $%
\mathbb{P}\{|X-m\mu|\geq c_{1}m\}$ is exponentially small [cf. Koml\'{o}s,
Major and Tusn\'{a}dy (\citeyear{KMT75})], an application of Lemma~\ref{KMT.lem} implies
$\mathbb{E}|\xi_{2}|^{k}\leq C_{k}m^{-k}$. Note that on the event $%
\{|X-m\mu|\leq c_{1}m\}$, $G^{\prime\prime} ( \mu^{\ast} ) $
is bounded for $m$ sufficiently large, then $\mathbb{E}|\xi
_{3}|^{k}\leq
C_{k}m^{-k}$ by observing that $\mathbb{E} [ ( X-m\mu) /%
\sqrt{m} ] ^{2k}\leq C_{k}^{\prime}$. The inequality $\mathbb{E}|\xi
|^{k}\leq C_{k}m^{-k}$ then follows immediately by combining the moments
bounds for $\xi_{1}$, $\xi_{2}$ and $\xi_{3}$. The second bound in
(\ref%
{moment.bd}) is a direct consequence of the first one and Markov inequality.
\end{pf*}

The variance stabilizing transformation considered in Section \ref%
{vst.sec} is for i.i.d. observations. In the function estimation procedure,
observations in each bin are independent but not identically distributed.
However, observations in each bin can be treated as i.i.d. random variables
through coupling. Let $X_{i}\sim \operatorname{NQ}(\mu_{i})$, $i=1,\ldots,m$, be independent.
Here the means $\mu_{i}$ are ``close'' but
not equal. Let $X_{i,c}$ be a set of i.i.d. random variables with
$X_{i,c}$ $%
\sim \operatorname{NQ}(\mu_{c})$. We define
\[
D=G \biggl( \frac{\sum_{i=1}^{m}X_{i}+a}{m+b} \biggr) -G \biggl( \frac{%
\sum_{i=1}^{m}X_{i,c}+a}{m+b} \biggr).
\]
If $\mu_{c}=\max_{i}\mu_{i}$, it is easy to see $\mathbb{E}D\leq0$
since $%
X_{i,c}$ is stochastically larger than $X_{i}$ for all $i$ [see, e.g.,
Lehmann and Romano (\citeyear{LR05})]. Similarly, $\mathbb{E}D\geq0$ when $\mu
_{c}=\min_{i}\mu_{i}$. We will select a
%
\begin{equation}\label{interthm}
\mu_{c}^{\ast}\in\Bigl[ \min_{i}\mu_{i},\max_{i}\mu_{i} \Bigr]
\end{equation}
such that $\mathbb{E}D=0$, which is possible by the intermediate value
theorem. In the following lemma we construct i.i.d. random variables $X_{i,c}
$ $\sim \operatorname{NQ}(\mu_{c}^{\ast})$ on the sample space of $X_{i}$ such that $D$
is very small and has negligible contribution to the final risk bounds in
Theorems \ref{global.thm} and \ref{LocalAdapt.thm}.
\begin{lemma}
\label{iidapprox}
Let $X_{i}\sim \operatorname{NQ}(\mu_{i})$, $i=1,\ldots,m$, be independent with
$\mu_{i}\in[ \varepsilon,v ]$,
a compact subset in the interior of the mean parameter space of
the natural exponential family. Assume that $ \vert
{\min_{i}\mu_{i}-\max_{i}\mu_{i}} \vert\leq C\delta$. Then
there are i.i.d. random variables $X_{i,c}$ where $X_{i,c}\sim
\operatorname{NQ}(\mu_{c}^{\ast})$ with $\mu_{c}^{\ast}\in[ \min_{i}\mu
_{i},\max_{i}\mu_{i} ] $ such that $\mathbb{E}D=0$ and:

\begin{longlist}
\item
%
\begin{equation}\label{tvapp}
\mathbb{P} ( \{ X_{i}\neq X_{i,c} \} ) \leq C\delta;
\end{equation}

\item and for any fixed integer $k\geq1$
there exists a constant $C_{k}>0$ such that for all $a>0$,
%
\begin{eqnarray} \label{moment.bd3}
\mathbb{E}|D|^{k} &\leq& C_{k}\log^{2k}m\cdot( m^{-k}+\delta
^{-k} ) \quad\mbox{and}\nonumber\\[-8pt]\\[-8pt]
\mathbb{P}(|D|>a) &\leq& C_{k}\frac{\log
^{2k}m}{a^{k}}(m^{-k}+\delta^{-k}).\nonumber
\end{eqnarray}
\end{longlist}
\end{lemma}
\begin{pf}
(i)
There is a classical coupling identity for the Total variation
distance. Let $P$ and $%
Q$ be distributions of two random variables $X$ and $Y$ on the same sample
space, respectively, then there is a random variable $Y_{c}$ with
distribution $Q$ such that $\mathbb{P} ( X\neq Y_{c} )
= \vert P-Q \vert_{\mathrm{TV}}$. See, for example, page 256 in Pollard
(\citeyear{Pollard02}). The proof of inequality (\ref{tvapp}) follows from that
identity and the inequality that $ \vert {\operatorname{NQ}(\mu_{i})-\operatorname{NQ}(\mu_{c}^{\ast
})} \vert_{\mathrm{TV}}\leq C \vert\mu_{i}-\mu_{c}^{\ast} \vert$
for some $C>0$ which only depends on the family of the distribution of $
X_{i} $ and $ [ \varepsilon,v ] $.

(ii) Using Taylor's expansion we can rewrite $D$ as $D=G^{\prime} (
\zeta) \frac{\sum_{i=1}^{m} ( X_{i}-X_{i,c} ) }{m+b}$ for
some $\zeta$ in between $\frac{\sum_{i=1}^{m}X_{i}+a}{m+b}$ and
$\frac{%
\sum_{i=1}^{m}X_{i,c}+a}{m+b}$. Since the distribution $X_{i}$ is in
exponential family, then $\mathbb{P} ( \max_{i} \vert
X_{i}-X_{i,c} \vert>\log^{2}m ) \leq C_{k^{\prime
}}m^{-k^{\prime}}$ for all $k^{\prime}>0$, which implies $\mathbb{E}
\vert X_{i}-X_{i,c} \vert^{k}\leq C_{k}\delta\log^{2k}m$ fo all
positive integer $k$. Since $X_{i}-X_{i,c}$ are independent, it can be shown
that
\begin{eqnarray*}
&&\mathbb{E} \Biggl( \frac{1}{m}\sum_{i=1}^{m} \vert
X_{i}-X_{i,c} \vert\Biggr) ^{k} \\
&&\qquad\leq\frac{1}{m^{k}}\sum_{k_{1}+\cdots+k_{m}=k}{\pmatrix{k\cr
k_{1},\ldots,k_{m}}%
}E \vert X_{1}-X_{1,c} \vert_{1}^{k_{1}}\cdots E \vert
X_{m}-X_{m,c} \vert_{m}^{k_{m}} \\
&&\qquad=\frac{1}{m^{k}}\sum_{j=1}^{k}\mathop{\sum_{k_{1}+\cdots+k_{m}=k,}}_
{\operatorname{Card} \{ i,k_{i}\geq1 \} =j}{\pmatrix{k\cr k_{1},\ldots,k_{m}}}
E \vert X_{1}-X_{1,c} \vert_{1}^{k_{1}}\cdots E \vert
X_{m}-X_{m,c} \vert_{m}^{k_{m}} \\
&&\qquad\leq C_{k}\frac{\log^{2k}m}{m^{k}}\sum_{j=1}^{k}\delta^{j}\cdot
\operatorname{Card} \bigl\{ ( k_{1},\ldots,k_{m} )\dvtx k_{1}+\cdots+k_{m}=k,
\\
&&\hspace*{200.4pt}
\operatorname{Card} \{ i,k_{i}\geq1 \} =j \bigr\} \\
&&\qquad\leq C_{k}^{\prime}\frac{\log^{2k}m}{m^{k}} \Biggl(
\sum_{j=1}^{k}m^{j}\delta^{j} \Biggr) =C_{k}^{\prime}\log^{2k}m \Biggl(
\sum_{j=1}^{k}m^{j-k}\delta^{j} \Biggr),
\end{eqnarray*}
where the last inequality follows from the facts that $k$ is fixed and
finite and
\begin{eqnarray*}
&&\operatorname{Card} \bigl\{ ( k_{1},\ldots,k_{m} ) \dvtx k_{1}+\cdots
+k_{m}=k, \operatorname{Card} \{ i,k_{i}\geq1 \} =j \bigr\} \\
&&\qquad = {\pmatrix{m\cr j}}\operatorname{Card} \{ ( k_{1},\ldots,k_{j} )
\dvtx k_{1}+\cdots+k_{j}=k, k_{i}\geq1 \} \\
&&\qquad \le{\pmatrix{m\cr j}} k^{k} \leq{m}^{j}k^{k}.
\end{eqnarray*}
Note that $\frac{m^{-k}+\delta^{k}}{m^{j-k}\delta^{j}}=\frac{1}{ (
m\delta) ^{j}}+ ( m\delta) ^{k-j}\geq1$ for all $k\geq
j\geq1$, then
\[
\mathbb{E} \Biggl( \frac{1}{m}\sum_{i=1}^{m} \vert
X_{i}-X_{i,c} \vert\Biggr) ^{k}\leq C_{k}^{\prime\prime}\log
^{2k}m\cdot( m^{-k}+\delta^{k} ).
\]
Thus the first inequality in (\ref{moment.bd3}) follows immediately by
observing that $G^{\prime} ( \zeta) $ is bounded with a
probability approaching to $1$ exponentially fast. The second bound is an
immediate consequence of the first one and Markov inequality.
\end{pf}
\begin{remark}
The unknown function $f$ in a Besov ball $B_{p,q}^{\alpha}( M) $ has
H\"{o}lder smoothness $d=\min(\alpha-{\frac{1}{p}},1)$, then $\delta
$ in Lemma \ref{iidapprox} can be chosen to be $T^{-d}$. The standard deviation
of normal noise in equation (\ref{Yi.decomp}) is $1/\sqrt{m}$. From the
assumptions in Theorems \ref{global.thm} or \ref{LocalAdapt.thm} we
see $%
m^{1/2}T^{-d}\log^{2}m$ converges to $0$ as a power of $n$, then
\begin{eqnarray*}
&&\mathbb{P}\bigl(|D|>1/\sqrt{m}\bigr)\\
&&\qquad\leq C_{k} \bigl[ ( m^{-1/2}\log^{2}m )
m^{-k}+ \bigl( \sqrt{m}T^{-d}\log^{2}m \bigr) ^{k} \bigr] \qquad\mbox{for all }%
k\geq1,
\end{eqnarray*}
which converges to $0$ faster than any polynomial of $m$. This implies the
contribution of $D$ to the final risk bounds in all major theorems is
negligible as shown in later sections.
\end{remark}

\subsection{\texorpdfstring{Proof of Theorem \protect\ref{Yi.thm}}{Proof of Theorem 1}}

From Lemma \ref{iidapprox}, there exist $Y_{j,c}^{\ast}$ where
$X_{i,c}\sim \operatorname{NQ}(f_{j}^{\ast})$ with
\[
f_{j,c}^{\ast}\in\biggl[ \min_{jm+1\leq i\leq( j+1 ) m}f \biggl(
\frac{i}{n} \biggr) ,\max_{jm+1\leq i\leq( j+1 ) m}f\biggl ( \frac{i%
}{n} \biggr) \biggr]
\]
as in (\ref{interthm}) such that%
%
\begin{eqnarray}
\label{err1}
\mathbb{E} [ Y_{j}^{\ast}-Y_{j,c}^{\ast} ] &=& 0, \\
\label{err21}
\mathbb{E}|Y_{j}^{\ast}-Y_{j,c}^{\ast}|^{k} &\leq& C_{k}\log
^{2k}m\cdot
( m^{-k}+T^{-dk} ), \\
\label{err2}
\mathbb{P}(|Y_{j}^{\ast}-Y_{j,c}^{\ast}|>a) &\leq& C_{k}\frac{\log
^{2k}m}{%
a^{k}} ( m^{-k}+T^{-dk} ).
\end{eqnarray}
Lemmas \ref{approx.lem}, \ref{KMT.lem} and \ref{xi.lem} together yield
%
\begin{equation} \label{err3}
Y_{j,c}^{\ast}=G(f_{j,c}^{\ast})+\epsilon_{j}+m^{-{
{1/2}}}Z_{j}+\xi
_{j},\qquad j=1,2,\ldots,T,
\end{equation}
and
%
\begin{equation} \label{err5}\qquad
\vert\epsilon_{j} \vert\leq Cm^{-2},\qquad \mathbb{E}|\xi
_{j}|^{k}\leq C_{k}m^{-k}\quad\mbox{and}\quad\mathbb{P}(|\xi_{j}|>a)\leq
C_{k} ( am ) ^{-k}.
\end{equation}
Note that%
%
\begin{equation} \label{err4}
\biggl\vert G(f_{j,c}^{\ast})-G \biggl( f\biggl(\frac{j}{T}\biggr) \biggr) \biggr\vert\leq
CT^{-d}.
\end{equation}
Theorem \ref{Yi.thm} then follows immediately by combining equations
(\ref%
{err1})--(\ref{err4}). 

\subsection{Risk bound for wavelet thresholding}
\label{SingleBlockRisk.sec}

We collect here a few technical results that are useful for the proof
of the
main theorems. We begin with the following moment bounds for an orthogonal
transform of independent variables. See Brown et al.
(\citeyear{Brownetal08}) for a proof.
\begin{lemma}
\label{moment.bound.lem} Let $X_{1},\ldots,X_{n}$ be independent variables
with $\mathbb{E}(X_{i})=0$ for $i=1,\ldots,n$. Suppose that $\mathbb
{E}%
|X_{i}|^{k}<M_{k}$ for all $i$ and all $k>0$ with $M_{k}>0$ some constant
not depending on $n$. Let $Y=WX$ be an orthogonal transform of $%
X=(X_{1},\ldots,X_{n})^{\prime}$. Then there exist constants
$M_{k}^{\prime}$
not depending on $n$ such that $\mathbb{E}|Y_{i}|^{k}<M_{k}^{\prime}$ for
all $i=1,\ldots,n$ and all $k>0$.
\end{lemma}

Lemma \ref{block.oracle.lem} below provides an oracle inequality for block
thresholding estimators without the normality assumption.
\begin{lemma}
\label{block.oracle.lem} Suppose $y_{i}=\theta_{i}+z_{i}, i=1,\ldots,L$,
where $\theta_{i}$ are constants and $z_{i}$ are random variables. Let
$%
S^{2}=\sum_{i=1}^{L}y_{i}^{2}$ and let $\widehat{\theta}_{i}=(1-{\frac
{\lambda L%
}{S^{2}}})_{+}y_{i}$. Then
%
\begin{equation} \label{block.oracle}
\mathbb{E}\Vert\widehat{\theta}-\theta\Vert_{2}^{2}\leq\Vert\theta
\Vert
_{2}^{2}\wedge4\lambda L+4\mathbb{E} [ \Vert z\Vert_{2}^{2}I(\Vert
z\Vert_{2}^{2}>\lambda L) ].
\end{equation}
\end{lemma}
\begin{pf} It is easy to verify that $%
\Vert\widehat{\theta}-y\Vert_{2}^{2}\leq\lambda L$. Hence
%
\begin{eqnarray} \label{tail.part}
&&\mathbb{E} [ \Vert\widehat{\theta}-\theta\Vert_{2}^{2}I(\Vert z\Vert
_{2}^{2}>\lambda L) ] \nonumber\\
&&\qquad\leq 2\mathbb{E} [ \Vert\widehat{\theta}%
-y\Vert_{2}^{2}I(\Vert z\Vert_{2}^{2}>\lambda L) ] +2\mathbb{E} [
\Vert y-\theta\Vert_{2}^{2}I(\Vert z\Vert_{2}^{2}>\lambda L) ]
\nonumber\\[-8pt]\\[-8pt]
&&\qquad\leq 2\lambda L\mathbb{P}(\Vert z\Vert_{2}^{2}>\lambda L)+2\mathbb
{E}%
[ \Vert z\Vert_{2}^{2}I(\Vert z\Vert_{2}^{2}>\lambda L) ]
\nonumber\\
&&\qquad\leq 4\mathbb{E} [ \Vert z\Vert_{2}^{2}I(\Vert z\Vert
_{2}^{2}>\lambda L) ].\nonumber
\end{eqnarray}
On the other hand,
%
\begin{eqnarray} \label{main.part1}
&&\mathbb{E} [ \Vert\widehat{\theta}-\theta\Vert_{2}^{2}I(\Vert z\Vert
_{2}^{2}\leq\lambda L) ] \nonumber\\[-8pt]\\[-8pt]
&&\qquad\leq\mathbb{E} [ (2\Vert\widehat{\theta}%
-y\Vert_{2}^{2}+2\Vert y-\theta\Vert_{2}^{2})I(\Vert z\Vert
_{2}^{2}\leq
\lambda L) ] \leq4\lambda L.\nonumber
\end{eqnarray}
Note that when $S^{2}\leq\lambda L$, $\widehat{\theta}=0$ and hence
$\Vert\widehat{%
\theta}-\theta\Vert_{2}^{2}=\Vert\theta\Vert_{2}^{2}$. When
$\Vert
z\Vert_{2}^{2}\leq\lambda L$ and $S^{2}>\lambda L$,
\begin{eqnarray*}
\Vert\widehat{\theta}-\theta\Vert_{2}^{2} &=&\sum_{i}\biggl[\biggl(1-{\frac
{\lambda L}{%
S^{2}}}\biggr)y_{i}-\theta_{i}\biggr]^{2}=\biggl(1-{\frac{\lambda
L}{S^{2}}}\biggr)\biggl[S^{2}-\lambda
L-2\sum_{i}\theta_{i}y_{i}\biggr]+\Vert\theta\Vert_{2}^{2} \\
&=&\biggl(1-{\frac{\lambda L}{S^{2}}}\biggr)\biggl[\sum(\theta_{i}+z_{i})^{2}-\lambda
L-2\sum_{i}\theta_{i}(\theta_{i}+z_{i})\biggr]+\Vert\theta\Vert_{2}^{2}
\\
&=&\biggl(1-{\frac{\lambda L}{S^{2}}}\biggr)(\Vert z\Vert_{2}^{2}-\lambda L-\Vert
\theta\Vert_{2}^{2})+\Vert\theta\Vert_{2}^{2}\leq\Vert\theta
\Vert
_{2}^{2}.
\end{eqnarray*}
Hence $\mathbb{E} [ \Vert\widehat{\theta}-\theta\Vert_{2}^{2}I(\Vert
z\Vert_{2}^{2}\leq\lambda L) ] \leq\Vert\theta\Vert_{2}^{2}$ and (
\ref{block.oracle}) follows by combining this with (\ref{tail.part})
and (%
\ref{main.part1}).
\end{pf}

The following bounds concerning a central chi-square distribution are from
Cai (\citeyear{Cai02}).
\begin{lemma}
\label{chisq.lem} Let $X\sim\chi_{L}^{2}$ and $\lambda>1$. Then
%
\begin{eqnarray} \label{chisq-tail}
\mathbb{P}(X\geq\lambda L)&\leq& e^{-{{L/2}}(\lambda-\log
\lambda
-1)} \quad\mbox{and}\nonumber\\[-8pt]\\[-8pt]
\mathbb{E}XI(X\geq\lambda L)&\leq&\lambda Le^{-{%
{L/2}}(\lambda-\log\lambda-1)}.\nonumber
\end{eqnarray}
\end{lemma}

From\vspace*{-2pt} (\ref{Yi.decomp}) in Theorem \ref{Yi.thm} we can write ${\frac
{1}{\sqrt{%
T}}}Y_{i}^{\ast}={\frac{G(f({i/T}))}{\sqrt{T}}}+{\frac{\epsilon
_{i}}{\sqrt{%
T}}}+{\frac{Z_{i}}{\sqrt{n}}}+{\frac{\xi_{i}}{\sqrt{T}}}$. Let $%
(u_{j,k})=T^{-{{1/2}}}W\cdot Y^{\ast}$ be the discrete wavelet
transform of the binned and transformed data. Then one may write
%
\begin{equation} \label{yjk.decomp}
u_{j,k}=\theta_{j,k}^{\prime}+\epsilon_{j,k}+{\frac{1}{\sqrt{n}}}%
z_{j,k}+\xi_{j,k},
\end{equation}
where $\theta_{jk}^{\prime}$ are the discrete wavelet transform of
$(G(f({%
i/T}))/\sqrt{T})$ which are approximately equal to the true wavelet
coefficients of $G ( f ) $, $z_{j,k}$ are the transform of the $%
Z_{i}$'s and so are i.i.d. $N(0,1)$ and $\epsilon_{j,k}$ and $\xi_{j,k}$
are, respectively, the transforms of $({\frac{\epsilon_{i}}{\sqrt
{T}}})$ and $%
({\frac{\xi_{i}}{\sqrt{T}}})$. Then it follows from Theorem \ref{Yi.thm}
that
%
\begin{equation} \label{epsilon}
\sum_{j}\sum_{k}\epsilon_{j,k}^{2}={\frac{1}{T}}\sum_{i}\epsilon
_{i}^{2}\leq C ( m^{-4}+T^{-2d} )
\end{equation}
and for all $i>0$ and $a>0$ we have%
%
\begin{eqnarray} \label{moment.bd2}
\mathbb{E}|\xi_{j,k}|^{i} &\leq& C_{i}^{\prime}\log^{2k}m \bigl[ (mn)^{-{%
{i/2}}}+T^{- ( d+1/2 ) i} \bigr],\nonumber\\[-8pt]\\[-8pt]
\mathbb{P}(|\xi_{j,k}|>a) &\leq& C_{i}^{\prime}\log^{2k}m [
(a^{2}mn)^{-{{i/2}}}+ ( aT^{d+1/2} ) ^{-i} ] \nonumber
\end{eqnarray}
from Theorem \ref{Yi.thm} and Lemma \ref{moment.bound.lem}.

Lemmas \ref{block.oracle.lem} and \ref{chisq.lem} together yield the
following result on the risk bound for a single block.
\begin{proposition}
\label{single.block.prop} Let the\vspace*{-1pt} empirical wavelet coefficients $%
u_{j,k}=\theta_{j,k}^{\prime}+\epsilon_{j,k}+{\frac{1}{\sqrt{n}}}%
z_{j,k}+\xi_{j,k}$ be given as in (\ref{yjk.decomp}) and let the block
thresholding estimator $\widehat{\theta}_{j,k}$ be defined as in (\ref
{block.est}%
). Then:

\begin{longlist}
\item for some constant $C>0$,
%
\begin{eqnarray}
\mathbb{E}\sum_{(j,k)\in B_{j}^{i}}(\widehat{\theta}_{j,k}-\theta
_{j,k}^{\prime
})^{2}&\leq&\min\biggl\{4\sum_{(j,k)\in B_{j}^{i}}(\theta_{j,k}^{\prime
})^{2}, 8\lambda_{\ast}Ln^{-1}\biggr\}\nonumber\\[-8pt]\\[-8pt]
&&{}+6\sum_{(j,k)\in B_{j}^{i}}\epsilon
_{j,k}^{2}+CLn^{-2};\nonumber
\end{eqnarray}

\item for any $0<\tau<1$, there exists a constant $C_{\tau}>0$
depending on $\tau$ only such that for all $(j,k)\in B_{j}^{i}$,
%
\begin{equation}\label{single.term.bnd}\quad
\mathbb{E}(\widehat{\theta}_{j,k}-\theta_{j,k}^{\prime})^{2}\leq
C_{\tau
}\cdot\min\Bigl\{ \max_{(j,k)\in B_{j}^{i}}\{(\theta_{j,k}^{\prime
}+\epsilon_{j,k})^{2}\}, Ln^{-1} \Bigr\} +n^{-2+\tau}.
\end{equation}
\end{longlist}
\end{proposition}

The following is a standard bound for wavelet approximation error. It
follows directly from Lemma 1 in Cai (\citeyear{Cai02}).
\begin{lemma}
\label{besov.approx.lem} Let $T=2^{J}$ and $d=\min(\alpha-{\frac
{1}{p}},1)$%
. Set
\[
\bar{g}_{J}(x)=\sum_{k=1}^{T}{\frac{1}{\sqrt{T}}}G \bigl( f (
k/n ) \bigr) \phi_{J,k}(x).
\]
Then for some constant $C>0$
%
\begin{equation}
{\sup_{g\in F_{p,q}^{\alpha}(M,\varepsilon)}}\Vert\bar{g}_{J}-G (
f ) \Vert_{2}^{2}\leq CT^{-2d}.
\end{equation}
\end{lemma}

We are now ready to prove our main results, Theorems \ref{global.thm}
and %
\ref{NEF-global.thm}.

\subsection{\texorpdfstring{Proofs of Theorems \protect\ref{global.thm} and \protect
\ref{NEF-global.thm}}{Proofs of Theorems 2 and 5}}

\label{hellinger.proof}

We shall only prove the results for the estimator $\widehat{f}_{\mathrm{BJS}}$. The proof
for $\widehat{f}_{\mathrm{NC}}$ is similar and simpler. Let $\widetilde{G ( f ) }%
=\max\{ \widehat{G ( f ) },0 \} $ for negative Binomial
and NEF--GHS distributions and $\widetilde{G ( f ) }=\widehat{%
G ( f ) }$ for other four distributions. We have%
\begin{eqnarray*}
\mathbb{E}\Vert\widehat{f}-f\Vert_{2}^{2} &=&\mathbb{E}\Vert
G^{-1}[%
\widetilde{G ( f ) }]-G^{-1}[G ( f ) ]\Vert_{2}^{2}=%
\mathbb{E}\Vert(G^{-1})^{\prime} ( g ) [\widetilde{G (
f ) }-G ( f ) ]\Vert_{2}^{2} \\
&\leq&\mathbb{E}\int V ( G^{-1} ( g ) ) [\widehat{%
G ( f ) }-G ( f ) ]^{2}\,dt,
\end{eqnarray*}
where $g$ is a function in between $\widetilde{G ( f ) }$ and $%
G ( f ) $. We will first give a lemma which implies $V (
G^{-1} ( g ) ) $ is bounded with high probability, then prove
Theorems~\ref{global.thm} and~\ref{NEF-global.thm} by establishing a risk
bound for estimating $G ( f ) $.
\begin{lemma}
\label{boundedest} Let $\widehat{G ( f ) }$ be the BlockJS estimator
of $G ( f ) $ defined in Section \ref{procedure.sec}. Then there
exists a constant $C>0$ such that
\[
\sup_{f\in F_{p,q}^{\alpha}(M,\varepsilon,v)}\mathbb{P} \{ \Vert
\widehat{G ( f ) } \Vert_{\infty}>C \} \leq C_{l}n^{-l}
\]
for any $l>1$, where $C_{l}$ is a constant depending on $l$.
\end{lemma}
\begin{pf}
Recall that we can write the discrete wavelet transform of the binned data
as
\[
u_{j,k}=\theta_{j,k}^{\prime}+\epsilon_{j,k}+{\frac{1}{\sqrt{n}}}%
z_{j,k}+\xi_{j,k},
\]
where $\theta_{jk}^{\prime}$ are the discrete wavelet transform of $({
\frac{G ( f ( i/T ) ) }{\sqrt{T}}})$ which are
approximately equal to the true wavelet coefficients $\theta_{jk}$ of $
G(f) $. Note that $ \vert\theta_{jk}^{\prime}-\theta_{jk} \vert
=O ( 2^{-j ( d+1/2 ) } ) ,\mbox{ for }d=\min( \alpha
-1/p,1 )$. Note also that a Besov Ball $B_{p,q}^{\alpha} (
M ) $ can be embedded in $B_{\infty,\infty}^{d} ( M_{1} ) $
for some $M_{1}>0$ [see, e.g., Meyer (\citeyear{meyer92})]. From the equation above, we
have
\[
\sum_{k=1}^{2^{j_{0}}}\widetilde{\theta}_{j_{0},k}^{\prime
}\phi_{j_{0},k}(t)+\sum_{j=j_{0}}^{J-1}\sum_{k=1}^{2^{j}}\theta
_{j,k}^{\prime}\psi_{j,k}(t)\in B_{\infty,\infty}^{d} ( M_{2} )
\]
for some $M_{2}>0$. Applying the Block thresholding approach, we have
\begin{eqnarray*}
\widehat{\theta}_{jk} &=&\biggl(1-{\frac{\lambda L\sigma^{2}}{S_{(j,i)}^{2}}}%
\biggr)_{+}\theta_{j,k}^{\prime}+\biggl(1-{\frac{\lambda L\sigma
^{2}}{S_{(j,i)}^{2}}}%
\biggr)_{+}\epsilon_{j,k}\\
&&{} +\biggl(1-{\frac{\lambda L\sigma^{2}}{S_{(j,i)}^{2}}}%
\biggr)_{+} \biggl( {\frac{1}{\sqrt{n}}}z_{j,k}+\xi_{j,k} \biggr) \\
&=&\widehat{\theta}_{1,jk}+\widehat{\theta}_{2,jk}+\widehat{\theta
}_{3,jk}\qquad\mbox{for }%
( j,k ) \in B_{j}^{i},j_{0}\leq j<J.
\end{eqnarray*}
Note that $ \vert\widehat{\theta}_{1,jk} \vert\leq\vert\theta
_{j,k}^{\prime} \vert$ and so $ \widehat{g}%
_{1}=\sum_{k=1}^{2^{j_{0}}}\widetilde{\theta}_{j_{0},k}^{\prime
}\phi
_{j_{0},k}+\sum_{j=j_{0}}^{J-1}\sum_{k=1}^{2^{j}}\widehat{\theta
}_{1,j,k}\psi
_{j,k}\in B_{\infty,\infty}^{d} ( M_{2} )$. This implies $%
\widehat{g}_{1}$ is uniformly bounded. Note that
\[
T^{{1/2}} \biggl(
\sum_{j,k} ( \epsilon_{j,k}^{2} ) \biggr) ^{1/2}=T^{{1/2}%
}\cdot O ( m^{-2} ) =o ( 1 ),
\]
so $W^{-1}\cdot T^{%
{1/2}} ( \widehat{\theta}_{2,jk} ) $ is a uniformly bounded
vector. For $0<\beta<1/6$ and a constant $a>0$ we have
\begin{eqnarray*}
\mathbb{P} \bigl( \vert\widehat{\theta}_{3,jk} \vert>a2^{-j (
\beta+1/2 ) } \bigr) &\leq&\mathbb{P} \bigl( \vert\widehat{\theta}%
_{3,jk} \vert>aT^{- ( \beta+1/2 ) } \bigr) \\
&\leq&\mathbb{P} \biggl( \biggl\vert{\frac{1}{\sqrt{n}}}z_{j,k} \biggr\vert>%
\frac{1}{2}aT^{- ( \beta+1/2 ) } \biggr) \\
&&{} +\mathbb{P} \biggl(
\vert\xi_{j,k} \vert>\frac{1}{2}aT^{- ( \beta+1/2 )
} \biggr) \leq A_{l}n^{-l}
\end{eqnarray*}
for any $l>1$ by Mill's ratio inequality and equation (\ref{moment.bd2}).
Let
\[
A=\bigcup_{j,k} \bigl\{ \vert\widehat{\theta}%
_{3,jk} \vert>a2^{-j ( \beta+1/2 ) } \bigr\}.
\]
Then $%
\mathbb{P} ( A ) =C_{l}n^{-l}$. On the event $A^{c}$ we have
\[
\widehat{g}_{3} ( t ) =\sum_{j=j_{0}}^{J-1}\sum_{k=1}^{2^{j}}\widehat{%
\theta}_{3,jk}\psi_{j,k}(t)\in B_{\infty,\infty}^{\beta} (
M_{3} ) \qquad\mbox{for some }M_{3}>0,
\]
which is uniformly bounded. Combining these results, we know that for $C$
sufficiently large
%
\begin{equation}\label{fbound}
\sup_{f\in F_{p,q}^{\alpha}(M,\varepsilon,v)}\mathbb{P} \{ \Vert
\widehat{G ( f ) } \Vert_{\infty}>C \} \leq\sup_{f\in
F_{p,q}^{\alpha}(M,\varepsilon)}\mathbb{P} ( A ) =C_{l}n^{-l}.
\end{equation}
\upqed\end{pf}

Now we are ready to prove Theorems \ref{global.thm} and \ref
{NEF-global.thm}%
. Note that $G^{-1}$ is an increasing and nonnegative function, and $V$
is a
quadratic variance function [see (\ref{Q.var})]. Lem\-ma~\ref%
{boundedest} implies that there exists a constant $C$ such that%
\[
\sup_{f\in F_{p,q}^{\alpha}(M,\varepsilon,v)}\mathbb{P} \{ \Vert
V ( G^{-1} ( g ) ) \Vert_{\infty}>C \} \leq
C_{l}n^{-l}
\]
for any $l>1$. Thus it is enough to show $\sup_{f\in F_{p,q}^{\alpha
}(M,\varepsilon,v)}\mathbb{E}\Vert\widehat{G ( f ) }-G (
f ) \Vert_{2}^{2}\leq Cn^{-{({2\alpha})/({1+2\alpha})}}$ for $p\geq
2 $ and $Cn^{-{({2\alpha})/({1+2\alpha})}}(\log n)^{
({2-p})/({p(1+2\alpha
)})}$ for $1\leq p<2$ under assumptions in Theorems \ref{global.thm}
and \ref%
{NEF-global.thm}.
\begin{pf*}{Proof of Theorem \protect\ref{global.thm}}
Let $Y$ and $\widehat{\theta}$ be given as in (\ref{data}) and
(\ref%
{block.est}), respectively. Then
%
\begin{eqnarray} \label{decomp3}
\mathbb{E}\Vert\widehat{G ( f ) }-G ( f ) \Vert_{2}^{2}
&=&\sum_{k}\mathbb{E}(\widehat{\!\widetilde{\theta}}_{j_{0},k}-{\widetilde
{\theta}}%
_{j,k})^{2}\nonumber\\
&&{} +\sum_{j=j_{0}}^{J-1}\sum_{k}\mathbb{E}(\widehat{\theta
}_{j,k}-\theta
_{j,k})^{2}+\sum_{j=J}^{\infty}\sum_{k}\theta_{j,k}^{2} \\
&\equiv&S_{1}+S_{2}+S_{3}.\nonumber
\end{eqnarray}
It is easy to see that the first term $S_{1}$ and the third term
$S_{3}$ are small:
%
\begin{equation}\label{s1}
S_{1}=2^{j_{0}}n^{-1}\epsilon^{2}=o\bigl(n^{-2\alpha/(1+2\alpha)}\bigr).
\end{equation}
Note that for $x\in\mathbb{R}^{m}$ and $0<p_{1}\leq p_{2}\leq\infty$,
%
\begin{equation} \label{2norms}
\Vert x\Vert_{p_{2}}\leq\Vert x\Vert_{p_{1}}\leq m^{{
{1}/{p_{1}}}-{%
{1}/{p_{2}}}}\Vert x\Vert_{p_{2}}.
\end{equation}
Since $f\in B_{p,q}^{\alpha}(M)$, so $2^{js}(\sum
_{k=1}^{2^{j}}|\theta
_{jk}|^{p})^{1/p}\leq M$. Now (\ref{2norms}) yields that
%
\begin{equation} \label{s3}
S_{3}=\sum_{j=J}^{\infty}\sum_{k}\theta_{j,k}^{2}\leq
C2^{-2J(\alpha
\wedge(\alpha+{{1/2}}-{{1/p}}))}.
\end{equation}
Proposition \ref{single.block.prop}, Lemma \ref{besov.approx.lem} and
(\ref{epsilon}) yield that
%
\begin{eqnarray} \label{S2}\hspace*{24pt}
S_{2} &\leq&2\sum_{j=j_{0}}^{J-1}\sum_{k}\mathbb{E}(\widehat{\theta}%
_{j,k}-\theta_{j,k}^{\prime})^{2}+2\sum_{j=j_{0}}^{J-1}\sum
_{k}(\theta
_{j,k}^{\prime}-\theta_{j,k})^{2} \nonumber\\
&\leq&\sum_{j=j_{0}}^{J-1}\sum_{i=1}^{2^{j}/L}\min\biggl\{ 8\sum
_{(j,k)\in
B_{j}^{i}}\theta_{j,k}^{2}, 8\lambda_{\ast}Ln^{-1}
\biggr\}\nonumber\\[-8pt]\\[-8pt]
&&{}
+6\sum_{j=j_{0}}^{J-1}\sum_{k}\epsilon
_{j,k}^{2}+Cn^{-1}+10\sum_{j=j_{0}}^{J-1}\sum_{k}(\theta
_{j,k}^{\prime
}-\theta_{j,k})^{2} \nonumber\\
&\leq&\sum_{j=j_{0}}^{J-1}\sum_{i=1}^{2^{j}/L}\min\biggl\{ 8\sum
_{(j,k)\in
B_{j}^{i}}\theta_{j,k}^{2}, 8\lambda_{\ast}Ln^{-1} \biggr\}
+Cm^{-4}+Cn^{-1}+CT^{-2d},\nonumber
\end{eqnarray}
which we now divide into two cases. First consider the case $p\geq2$.
Let $J_{1}=[%
{\frac{1}{1+2\alpha}}\log_{2}n]$. So, $2^{J_{1}}\approx
n^{1/(1+2\alpha)}$%
. Then (\ref{S2}) and (\ref{2norms}) yield
%
\begin{eqnarray}\label{s2aa2}
S_{2}&\leq&8\lambda_{\ast
}\sum_{j=j_{0}}^{J_{1}-1}\sum_{i=1}^{2^{j}/L}Ln^{-1}+8\sum
_{j=J_{1}}^{J-1}%
\sum_{k}\theta_{j,k}^{2}+Cn^{-1}+CT^{-2d}\nonumber\\[-8pt]\\[-8pt]
&\leq& Cn^{-2\alpha
/(1+2\alpha)}.\nonumber
\end{eqnarray}
By combining (\ref{s2aa2}) with (\ref{s1}) and (\ref{s3}), we have
$\mathbb{E}%
\Vert\widehat{\theta}-\theta\Vert_{2}^{2}\leq Cn^{-2\alpha/(1+2\alpha)}$,
for $p\geq2$.
\end{pf*}

Now let us consider the case $p<2$. First we state the following lemma
without proof.
\begin{lemma}
\label{ubound0.lem} Let $0<p<1$ and $S=\{x\in\mathbb{R}^{k}\dvtx
\sum_{i=1}^{k}x_{i}^{p}\leq B, x_{i}\geq0, i=1,\ldots,k\}$. Then $%
\sup_{x\in S}\sum_{i=1}^{k}(x_{i}\wedge A)\leq B\cdot A^{1-p}$ for
all $A>0$.
\end{lemma}

Let $J_{2}$ be an integer satisfying $2^{J_{2}}\asymp
n^{1/(1+2\alpha)}(\log n)^{(2-p)/p(1+2\alpha)}$. Note that
\[
\sum_{i=1}^{2^{j}/L} \biggl( \sum_{(j,k)\in B_{j}^{i}}\theta_{j,k}^{2} \biggr)
^{{p/2}}\leq\sum_{k=1}^{2^{j}}(\theta_{j,k}^{2})^{
{p/2}}\leq
M2^{-jsp}.
\]
It then follows from Lemma \ref{ubound0.lem} that
%
\begin{eqnarray}\label{ss21}
&&\sum_{j=J_{2}}^{J-1}\sum_{i=1}^{2^{j}/L}\min\biggl\{ 8\sum_{(j,k)\in
B_{j}^{i}}\theta_{j,k}^{2}, 8\lambda_{\ast}Ln^{-1} \biggr\}\nonumber\\[-8pt]\\[-8pt]
&&\qquad \leq Cn^{-{%
({2\alpha})/({1+2\alpha})}}(\log
n)^{({2-p})/({p(1+2\alpha)})}.\nonumber
\end{eqnarray}
On the other hand,
%
\begin{eqnarray} \label{ss22}
&&\sum_{j=j_{0}}^{J_{2}-1}\sum_{i=1}^{2^{j}/L}\min\biggl\{ 8\sum_{(j,k)\in
B_{j}^{i}}\theta_{j,k}^{2}, 8\lambda_{\ast}Ln^{-1} \biggr\}\nonumber\\
&&\qquad\leq
\sum_{j=j_{0}}^{J_{2}-1}\sum_{b}8\lambda_{\ast}Ln^{-1}\\
&&\qquad\leq
Cn^{-{({%
2\alpha})/({1+2\alpha})}}(\log n)^{({2-p})/({p(1+2\alpha)})}.\nonumber
\end{eqnarray}
Putting (\ref{s1}), (\ref{s3}), (\ref{ss21}) and (\ref{ss22}) together
yields $\mathbb{E}\Vert\widehat{\theta}-\theta\Vert_{2}^{2}\leq
Cn^{-{({%
2\alpha})/({1+2\alpha})}}\times\break(\log n)^{({2-p})/({p(1+2\alpha)})}$.
%
%
\begin{pf*}{Proof of Theorem \protect\ref{NEF-global.thm}}
The proof of Theorem \ref{NEF-global.thm} is similar to that of
Theorem~\ref{global.thm} except the step of  (\ref{S2}). We will thus omit
most of the details. For a general natural exponential family the upper
bound for $\sum_{j=j_{0}}^{J-1}\sum_{k}\epsilon_{j,k}^{2}$ in
equation (\ref%
{S2}) is $C ( m^{-2}+T^{-2d} ) $ as given in Section \ref{vst.sec},
so (\ref{S2}) now becomes%
\[
S_2\le\sum_{j=j_{0}}^{J-1}\sum_{i=1}^{2^{j}/L}\min\biggl\{ 8\sum
_{(j,k)\in
B_{j}^{i}}\theta_{j,k}^{2}, 8\lambda_{\ast}Ln^{-1} \biggr\}
+Cm^{-2}+Cn^{-1}+CT^{-2d}.
\]
For $m=cn^{-1/2}$, we have $m^{-2}=c^{2}n^{-1}$. When $\alpha-{\frac
{1}{p}}>%
{\frac{2\alpha}{1+2\alpha}}$, it is easy to see $T^{-2d}=o (
n^{-2\alpha/(1+2\alpha)} ) $. Theorem \ref{NEF-global.thm} then
follows from the same steps as in the proof of Theorem \ref{global.thm}.
\end{pf*}

\section*{Acknowledgments}
We would like to thank Eric Kolaczyk for providing the BATSE data and Packet
loss data. We also thank the two referees for detailed and constructive
comments which have lead to significant improvement of the paper.

\printaddresses

\end{document}